\documentclass[10pt, a4paper]{article}
\usepackage{amssymb,amsmath,amsthm,tocloft,mathtools,comment}
\usepackage{fancyhdr}
\usepackage{enumerate}
\usepackage[british]{babel}
\usepackage{enumitem,verbatim}

\usepackage{stackengine}
\usepackage{hyperref}
\usepackage{cleveref}
\usepackage[texencoding=ascii,
    sorting=nyt,
    maxbibnames=99
]{biblatex}

\usepackage[usenames,dvipsnames,table]{xcolor}
\usepackage{graphicx,tikz,caption,subcaption}
\usetikzlibrary{decorations.pathreplacing}
\usetikzlibrary{decorations.pathmorphing}
\usetikzlibrary{calc,arrows,decorations.markings}
\usetikzlibrary{shapes}
\tikzset{snake it/.style={decorate, decoration=snake}}
\usepackage{bbm,dsfont}
\newcounter{propcounter}

\usepackage[margin=2.63cm]{geometry}

\newtheorem{theorem}{Theorem}[section]

\newtheorem{lemma}[theorem]{Lemma}

\newcommand{\repeatlabel}{}
\newtheorem*{repeatlemma}{Lemma \repeatlabel}

\theoremstyle{definition}

\newtheorem{defn}[theorem]{Definition}

\newtheorem{claim}[theorem]{Claim}

\newenvironment{poc}{\begin{proof}[Proof of claim]}{\end{proof}}

\newcommand{\eps}{\varepsilon}

\newcommand{\G}{\mathcal{G}}

\def \p {\phi}

\def \a {\alpha}
\def \g {\gamma}
\def \b {\beta}
\def \P {\mathbb{P}}
\def \E {\mathbb{E}}

\def \N {\mathbb{N}}

\def \M {\mathcal{M}}

\def \D {\Delta}
\def \d {\delta}
\def \e {\varepsilon}
\def \h {\eta}

\def \n {\tilde{n}}

\def \tV {\tilde{V}}

\title{A bandwidth theorem for graph transversals}

\author{Debsoumya Chakraborti\thanks{
Discrete Mathematics Group (DIMAG), Institute for Basic Science (IBS), Daejeon, South Korea.
E-mail: {\tt debsoumya@ibs.re.kr}. This work was supported by the Institute for Basic Science (IBS-R029-C1).}\and
Seonghyuk Im\thanks{Department of Mathematical Sciences, KAIST, South Korea and Extremal Combinatorics and Probability Group (ECOPRO), Institute for Basic Science (IBS), Daejeon, South Korea. E-mail: {\tt seonghyuk@kaist.ac.kr}. Seonghyuk Im was supported by the POSCO Science Fellowship of POSCO TJ Park Foundation and by the Institute for Basic Science (IBS-R029-C4).}\and
Jaehoon Kim\thanks{Department of Mathematical Sciences, KAIST, South Korea. E-mail: {\tt jaehoon.kim@kaist.ac.kr}. Jaehoon Kim was supported by the POSCO Science Fellowship of POSCO TJ Park Foundation.}\and
Hong Liu\thanks{Extremal Combinatorics and Probability Group (ECOPRO), Institute for Basic Science (IBS), Daejeon, South Korea. E-mail: {\tt hongliu@ibs.re.kr}. Supported by IBS-R029-C4.} }

\begin{document}

\maketitle

\begin{abstract}
Given a collection $\mathcal{G}=(G_1,\dots, G_h)$ of graphs on the same vertex set $V$ of size $n$, an $h$-edge graph $H$ on the vertex set $V$ is a $\mathcal{G}$-transversal if there exists a bijection $\lambda : E(H) \rightarrow [h]$ such that $e\in E(G_{\lambda(e)})$ for each $e\in E(H)$. The conditions on the minimum degree $\delta(\mathcal{G})=\min_{i\in[h]}\{ \delta(G_i)\}$ for finding a spanning $\mathcal{G}$-transversal isomorphic to a graph $H$ have been actively studied when $H$ is a Hamilton cycle, an $F$-factor, a spanning tree with maximum degree $o(n/\log n)$ and a power of a Hamilton cycle, etc. In this paper, we determined the asymptotically tight threshold on $\delta(\mathcal{G})$ for finding a $\mathcal{G}$-transversal isomorphic to $H$ when $H$ is a general $n$-vertex graph with bounded maximum degree and $o(n)$-bandwidth. This provides a transversal generalization of the celebrated Bandwidth theorem by B\"ottcher, Schacht and Taraz.
\end{abstract}


\section{Introduction}
For a given collection $\mathcal{F}=(F_1,\dots, F_s)$ of objects, a set $X$ which intersects with every $F_i$ is called a `transversal' of $\mathcal{F}$ or `colorful' set (with respect to $\mathcal{F}$).
Transversals of a collection of various mathematical objects have been extensively studied.
For example,  B\'ar\'any~\cite{Barany1982}  proved a colorful generalization of Carath\'eodory's theorem and Holmsen, Pach and Tverberg \cite{Holmsen2008} further generalized it. \v{Z}ivaljevi\'{c} and  Vre\'{c}ica~\cite{Zivaljevic1992} proved colorful version of Tverberg's theorem. Kalai and Meshulam~\cite{Kalai2005} considered a topological colorful version of Helly's theorem. Aharoni and Howard \cite{Aharoni2017} proved a colorful version of the Erd\H{o}s-Ko-Rado theorem.
Rota \cite{Huang1994} conjectured that a collection of $n$ bases of an $n$-dimensional vector space can be decomposed into $n$ transversals so that each of which is also a basis. This was further generalized into Rota's basis conjecture stating that any collection of $n$ bases of rank $n$ matroids can be decomposed into transversals so that each of them is also a base. There have been vibrant progress on this conjecture; for example, an asymptotic version of this conjecture was recently solved by Pokrovskiy in \cite{Pokrovskiy2020}.

Even though it is a very natural concept, the transversals of collections of graphs were not formally introduced until the recent work of Joos and the third author \cite{Joos2020}. It reads as follows.
\begin{defn}
For a given collection $\mathcal{G}=(G_1,\dots, G_h)$ of graphs with the same vertex set $V$, an $h$-edge graph $F$ on the vertex set $V$ is a $\mathcal{G}$-transversal if there exists a bijection $\lambda:E(F)\rightarrow [h]$ such that $e\in E(G_{\lambda(e)})$ for all $e\in E(F)$.
\end{defn}

By considering each $G_i$ as a graph whose edges are colored with $i$, the function $\lambda$ is often called a \emph{coloring}. We call this coloring \emph{rainbow} when $\lambda$ is injective. We call a graph equipped with a rainbow coloring a \emph{partial transversal}. 
Even before its formal definition, the (partial) transversals of graphs were widely studied, under the name of rainbow subgraph, for the cases when the graphs $G_1,\dots, G_h$ are edge-disjoint. In such cases, the graph collection $(G_1,\dots, G_h)$ is equivalent to an edge-colored graph, and one can investigate sufficient conditions to find a rainbow subgraph $H$ in edge-colored graphs.
For example, it was considered when the graphs $H$ are perfect matchings~\cite{Hatami2008,Pokrovskiy2018}, Hamilton cycles~\cite{Alon2017,Balogh2019,Benzing2020},  spanning trees~\cite{Balogh2018,Benzing2020,Sudakov2018} or other graphs~\cite{Bottcher2012,Ehard2020,Glock2020,Sudakov2017}.

Considering general cases where the graphs $G_1,\dots, G_h$ are not necessarily edge-disjoint also leads to very interesting phenomena, extending
many classical results in extremal graph theory. Aharoni, DeVos, de la Maza, Montejano and \v{S}\'{a}mal \cite{Aharoni2020} proved that if a graph collection $\mathcal{G}=(G_1,G_2,G_3)$ satisfies $e(G_i) > (\frac{26-2\sqrt{7}}{81})n^2$ for every $i\in [3]$, then it has a $\mathcal{G}$-transversal isomorphic to a triangle. Surprisingly, this condition is best possible showing an interesting behavior on $\mathcal{G}$-transversal as $\frac{26-2\sqrt{7}}{81}$ is larger than $\frac{1}{4}$ which we obtain from the Mantel's theorem. It is an interesting open problem to generalize this further by determining the tight conditions on $e(G_i)$ for the existence of a $\mathcal{G}$-transversal isomorphic to $K_r$ with $r>3$.
 
Motivated from a question in \cite{Aharoni2020}, 
Cheng, Wang and Zhao \cite{Cheng2021} proved an asymptotic extension of Dirac's theorem \cite{Dirac1952} and 
Joos and the third author \cite{Joos2020} proved the exact version by showing that if $\mathcal{G}=(G_1,\dots, G_n)$ satisfies $\delta(\mathcal{G})=\min_{i\in [n]} \{ \delta(G_i) \}\geq \frac{n}{2}$ contains a $\mathcal{G}$-transversal isomorphic to a Hamilton cycle. 
Montgomery, M\"{u}yesser and Pehova~\cite{Montgomery2022} determined asymptotically tight bound on $\delta(\mathcal{G})$ for the existence of a $\mathcal{G}$-transversal isomorphic to a general $F$-factor or a spanning tree with maximum degree at most $o(\frac{n}{\log n})$, generalizing the classical Hajnal-Szemer\'edi theorem~\cite{Hajnal1970} and the tree-embedding theorem of Koml\'os, S\'ark\"ozy and Szmer\'edi~\cite{Komlos2001} respectively. 
A similar result for $K_t$-factor was proved by
Cheng, Han, Wang and Wang 
\cite{Cheng2021a} and generalized into the cases of digraphs and hypergraphs, and 
Cheng, Wang and Zhao \cite{Cheng2021} further proved results for Hamilton cycles in the collections of hypergraphs.
Moreover, Gupta, Hamann, M\"uyesser, Parczyk and Sgueglia~\cite{Gupta2022} recently determined the asymptotically tight bound on $\delta(\mathcal{G})$ for the existence of powers of Hamilton cycles, generalizing P\'osa-Seymour conjecture, which was proved by Koml\'{o}s, S\'{a}rk\"{o}zy, and Szemer\'{e}di~\cite{Komlos1998}.
Cheng and Staden \cite{CS2023} developed a version of rainbow blow-up lemma (which can be used when the number of colors is $\eps$-fraction more than the number of edges in $H$) and obtained a result similar to \cite{Gupta2022} when the number of colors is $\eps$-fraction more than the number of edges in the power of Hamilton cycle. 

 All those graphs above, powers of Hamilton cycles, $F$-factors, trees, have somewhat bounded maximum degree and have low connectivity. This low connectivity can be captured by the following notion of bandwidth. A graph $H$ has a \textit{bandwidth} at most $b$ if there exists an ordering $x_1, \cdots, x_n$ of $V(H)$ such that all edges $x_ix_j \in E(H)$ satisfies $|i-j|\le b$.
Indeed, the celebrated 
 bandwidth theorem proved by B\"{o}ttcher, Schacht and Taraz \cite{Boettcher2009} determines the asymptotically sharp minimum degree condition on $G$ to find such a graph $H$ with bounded maximum degree and low bandwidth as a spanning subgraph.
More precisely, the bandwidth theorem states that if an $n$-vertex $k$-chromatic graph $H$ has bounded maximum degree and sublinear bandwidth, then every $n$-vertex graph $G$ with $\delta(G) \geq (1-\frac{1}{k}+o(1))n$ contains a copy of $H$. 

As the bandwidth theorem provides a common generalization for many classical results, it is very natural to pursue such a generalization for the above transversal-version results.
Indeed, we prove that such a generalization holds by proving the following bandwidth theorem for graph transversals. In short, it determines the asymptotically tight bound on $\delta(\mathcal{G})$ for finding a spanning $\mathcal{G}$-transversal isomorphic to $H$, when $H$ is any graph satisfying the conditions on the bandwidth theorem. 
Note that the original bandwidth theorem can be recovered by taking each $G_i$ to be the same graph. Thus, Theorem~\ref{main} generalizes and strengthens the bandwidth theorem.

\begin{theorem} \label{main}
For every $\eps >0$ and positive integers $\Delta, k$, there exist $\alpha>0$ and $h_0>0$ satisfying the following for every $h\ge h_0$. 
Let $H$ be an $n$-vertex graph with $h$ edges and bandwidth at most $\alpha n$ such that $\Delta(H)\leq \Delta$ and $\chi(H)\leq k$. If $\G = (G_1,\ldots,G_h)$ is a family of $h$ graphs on the same vertex set of size $n$ such that $\d(G_i)\geq  \left(1 - \frac{1}{k} + \eps\right)n$ for all $i\in [h]$, then there exists a $\G$-transversal isomorphic to~$H$.
\end{theorem}


\section{Preliminaries}

For $n\leq m\in \mathbb{N}$, we write $[n]:=\{1,\dots, n\}$ and $[n,m]:=\{ n, n+1, \dots, m\}$.
For $a,b,c\in \mathbb{R}$, we write $a= b\pm c$ to denote $b-c\leq a\leq b+c$.
We often treat large numbers as integers whenever it does not affect the argument.
If we claim that a statement holds when $0<1/n\ll a\ll b,c\leq 1$, then it means that the statement holds for $a,b,c, n$ satisfying $1/n <f(a)$ and $a< g(b,c)$ for some non-decreasing functions $f: (0,1]\rightarrow (0,1]$ and $g:(0,1]^2\rightarrow (0,1]$. We will not explicitly state these functions. Hierarchies with more constants are defined in a similar way. 
A set of size $k$ is called a $k$-set.
When a collection of sets $X_1,\dots, X_k$ and $Y\subseteq [k]$ are given, we write $X_Y = \bigcup_{i\in Y} X_i$ to denote the union of the sets $X_i$ for all $i\in Y$. 
For a set $X$, we denote by $\binom{X}{k}$ the family of $k$-element subsets of $X$.

We use standard graph theory notations. 
For a graph $H$, we denote the vertex set of $H$ by $V(H)$, the edge set by $E(H)$, the number of edges by $e(H)$, the maximum degree by $\D(H)$, the minimum degree by $\d(H)$ and the chromatic number by $\chi(H)$. For a graph $H$ and two disjoint subsets $A,B\subseteq V(H)$, we denote the induced subgraph on $A$ by $H[A]$, and the induced bipartite subgraph with parts $A,B$ by $H[A,B]$.

For a graph $G$, we write $K_k(G)$ to denote the set of all $k$-cliques in $G$ and we write
$\overrightarrow{K_k}(G)$ to denote the set of all ordered $k$-tuples $(v_1,\dots, v_k) \in V(G)^k$ where $G[\{v_1,\dots, v_k\}]$ forms a $k$-clique.
For a set $X\subseteq V(G)$, we write $N_G(X) = \bigcap_{x\in X} N_G(x)$ for the common neighborhood of the vertices and $N^*_G(X) = \bigcup_{x\in X}N_G(x)$.
We say that a partition $(V_1,\dots, V_r)$ of $V$ is \emph{equitable} if $||V_i|-|V_j||\leq 1$ for all $i,j\in [r]$.

Consider a graph $R$ on the vertex set $[r]$.
For a graph $H$ and a partition $\mathcal{X}=(X_1,\dots, X_r)$ of $V(H)$, we say that $H$ admits the vertex partition $(R,\mathcal{X})$ if each $X_i$ is an independent set in $H$ and  there are no edges in  $H[X_i,X_j]$ for each $ij\notin E(R)$. 

We will use the following concentration inequality. 

\begin{lemma}[Chernoff-Hoeffding inequality, Lemma~2.1 in \cite{Kwan2019}] \label{lem: sharp concentration}
Suppose that $g\colon \{0,1\}^n \rightarrow \mathbb{R}$ satisfies the bounded differences condition with parameters $(c_1,\ldots,c_n) \in \mathbb{R}^n$, i.e.,
\begin{align*} 
    |g(x)-g(x')|\le \sum_{i=1}^n c_i\mathbbm{1}_{\{x_i\neq x'_i\}} \quad \text{for every} \enskip x,x' \in \{0,1\}^n .
\end{align*}
Let $\xi \in \{0,1\}^n$ be a random vector uniformly distributed on $\binom{[n]}{m}$. Then, for each $t\ge 0$, 
$$\P\Big[\left|g(\mathbb{\xi}) - \E \big[g(\mathbb{\xi})\big]\right| \ge t\Big]\le 2e^{-\frac{t^2}{8\sum_{i=1}^n c_i^2}}.$$
\end{lemma}

Instead of using this lemma directly, it will be convenient to use the following simple corollary.
\begin{lemma} \label{lem:concentration}
Let $0<m,r\le n$. Suppose that $S \subseteq [n]$ is a fixed subset of size $m$ and that $\mathcal{X}$ is an $r$-set chosen uniformly at random from $\binom{[n]}{r}$. Then, for each $t\ge 0$, $$\P\left[\left||S \cap \mathcal{X}| - \frac{mr}{n}\right| \ge t\right] \le 2\exp\left(-\frac{t^2}{8r}\right).$$
\end{lemma}

\begin{proof}
Let $S,X\subseteq [n]$ be two fixed subsets of size $m$ and $r$ respectively. Let $\mathcal{S},\mathcal{X}$ be two sets chosen uniformly at random respectively from $\binom{[n]}{m}$ and $\binom{[n]}{r}$. Consider the random variables $\left|\mathcal{S} \cap \mathcal{X}\right|$, $\left|S \cap \mathcal{X}\right|$, and $\left|\mathcal{S} \cap X\right|$. It is easy to see that these three random variables have the same distribution. Thus, it is enough to prove that $$\P\left[\left||\mathcal{S} \cap X| - \frac{mr}{n}\right| \ge t\right] \le 2\exp\left(-\frac{t^2}{8r}\right).$$
This is a straightforward consequence of Lemma~\ref{lem: sharp concentration} applied with the parameters $(c_1,\dots, c_n)$ such that $c_i=1$ if $i\in X$ and $c_i=0$ otherwise.
\end{proof}

The minimum degree condition on a graph $H$ ensures that we can find a copy of $K_k$. For example, the following lemma is easy to prove.
\begin{lemma}\label{lem: clique exist}
If $H$ is an $n$-vertex graph with minimum degree at least $(1-\frac{1}{k-1})n+1$ and $1\leq k'\leq k$, then every copy of $K_{k'}$ in $H$ is contained in a copy of $K_k$.
\end{lemma}

The following Hajnal-Szemer\'edi theorem provides a $K_k$-factor in a graph with large minimum degree.
\begin{theorem}[Hajnal-Szemer\'edi Theorem \cite{Hajnal1970}]\label{thm:HS}
 If $G$ is an $n$-vertex graph with minimum degree at least $\left(1-\frac{1}{k}\right)n$ and $n$ is a multiple of $k$, then $G$ contains a $K_k$-factor.
\end{theorem}

We say that a matching $M$ in a graph $H$ is \emph{$3$-independent} if every two edges in $M$ are in the distance at least three apart. In other words, no vertex in $H$ has two neighbors $x,y$ covered by $M$, unless $xy\in E(M)$.
The following lemma yields a large $3$-independent matching in a graph $H$. 

\begin{lemma}\label{lem: matching M}
Let $H$ be a graph with maximum degree at most $\Delta$. 
Then any edge set $E$ of $H$ contains a $3$-independent matching of size at least $\frac{|E|}{2\Delta^3}$.
\end{lemma}

This can be proved by observing that such $3$-independent matchings are in correspondence with the independent sets in the $3$rd power of the line graph of $H$, which has maximum degree at most $2(\Delta -1)^3 + 2(\Delta-1)^2 + 2(\Delta-1) \le 2\Delta^3-1$. 


\subsection{Bandwidth ordering and clique-walks}
Consider a graph $H$ with bandwidth at most $b$. Let $(x_i)_{i\in [n]}= (x_1,\dots, x_n)$ be a $b$-bandwidth ordering of $V(H)$, i.e., if $x_ix_j\in E(H)$, then $|i-j|\leq b$.

If $H$ has small bandwidth, then the structure of $H$ resembles a path. This fact is reflected by the following lemma that provides us with a good partition of $V(H)$. Let $P_{[r]}^k$ be the $k$-th power of the path on the vertex set $[r]$ where $ij\in E(P_{[r]}^k)$ if and only if $|i-j|\leq k$.  We omit the proof of this lemma as it is  straightforward to check.
\begin{lemma}\label{lem: bandwidth partition}
Let $H$ be an $n$-vertex graph and $\chi(H)=k$ and it has an $\alpha n$-bandwidth ordering $(x_i)_{i\in [n]}$. Let $c:V(H)\rightarrow [k]$ be a proper $k$-coloring of $H$. 
Then the vertex partition $\mathcal{W}=(W_1,\dots, W_r)$ of $H$ with 
$$W_{i} = \{ x_{\ell}\in V(H) : (i-k)\alpha n +1 \leq \ell \leq i\alpha n, \; c(x_{\ell}) \equiv i ~({\rm mod}~ k)\}$$
satisfies the following.
\begin{itemize}
\item[(1)] For each $i\in [r]$, the set $W_i$ is an independent set with $|W_i|\leq k\alpha n$.
\item[(2)] $r$ is a multiple of $k$ and $H$ admits the vertex partition $\left(P_{[r]}^{k-1},\mathcal{W}\right)$. \end{itemize}
\end{lemma}
Note that some sets $W_i$ above can be empty.
In this lemma, the second condition means that we have $N^*(W_i) \subseteq \bigcup_{|j-i|\leq k-1} W_j$ for all $i\in [r]$. We call such a vertex partition $\mathcal{W}$ an \emph{$\alpha n$-bandwidth partition} of $H$ with respect to the bandwidth ordering $(x_i)_{i\in [n]}$ and the proper $k$-coloring $c$. Often when the ordering $(x_i)_{i\in [n]}$ and the coloring $c$ are understood, then we just say $\mathcal{W}$ is an $\alpha n$-bandwidth partition.

The following lemma yields a subgraph of $R$ which is roughly a $k$-th power of a walk. This lemma was implicitly used in \cite{Komlos1998} and was explicitly stated in \cite{Condon2019}.

\begin{lemma}[\cite{Condon2019, Komlos1998}]\label{lem: clique_walk} 
Let $r,k \in \N \setminus \{1\}$. Suppose that $R$ is an $r$-vertex graph with $\d(R) \ge \left(1-\frac{1}{k}\right)r+1$. Suppose that $Q_1,Q_2$ are two ordered tuples of $k$ distinct vertices such that $Q_1 = (x_1,\ldots,x_k)$ and $Q_2 = (y_1,\ldots,y_k)$. 
Two sets $\{x_1,\dots, x_k\}$ and $\{y_1,\dots, y_k\}$ do not have to be disjoint. Then there exists a walk $W = (z_1,\ldots,z_t)$ in $R$ satisfying the following.
\begin{itemize}
    \item $3k \le t \le 3k^3$ and $k$ divides $t$,
    \item for all $i,j \in [t]$ with $|i-j| \le k-1$, if $\{i,j\}\not\subseteq \{1,\dots, k, t-k+1,\dots,t\}$ we have $z_iz_j \in E(R)$,
    \item for each $i \in [k]$, we have $z_i=x_i$ and $z_{t-k+i}=y_i$.
\end{itemize}
\end{lemma}
If we have a walk $W$ as above, and if both $Q_1, Q_2$ form a clique in addition, then we call such a walk $W$ as a \emph{$k$-clique-walk} in $R$. We often just call it a \emph{clique-walk} if $k$ is well understood.

\subsection{{\boldmath $\varepsilon$}-regularity}
The notion of $\eps$-regularity is important in this paper.
We say that a bipartite graph $G$ with vertex partition $(A,B)$ is $(\eps,d)$-regular if all sets $A'\subseteq A$ and $B'\subseteq B$ with $|A'|\ge \eps|A|$ and $|B'|\ge \eps|B|$ satisfy 
$$\left| \frac{e_G(A',B')}{|A'||B'|} - d \right| \leq \eps.$$
If $G$ is $(\eps,d)$-regular for some $d\geq 0$, then it is called $\eps$-regular.
If $G$ is $(\eps,d')$-regular for some $d'\geq d$, then we say that it is $(\eps,d+)$-regular. If the induced bipartite graph $G[A,B]$ is $(\eps,d+)$-regular, then we say that $(A,B)$ is an $(\eps,d+)$-regular pair in $G$.
We say that $(A,B)$ is an $(\eps,d)$-super-regular pair in $G$, if it is an $(\eps,d+)$-regular pair and additionally, all vertices $v\in A$ and $u\in B$ satisfy $|N_G(v)\cap B|\geq (d-\eps)|B|$ and $|N_G(u)\cap A|\geq (d-\eps)|A|$.

\begin{defn}
Let $R$ be a graph on the vertex set $[r]$ and $\mathcal{V}=(V_1,\dots, V_r)$ be a partition of $V(G)$.
We say that $G$ admits \emph{an $(\eps,d)$-regularity partition $(R,\mathcal{V})$}	if $G[V_i,V_j]$ is $(\eps,d+)$-regular for all $ij\in E(R)$. 
We say that $G$ admits \emph{an $(\eps,d)$-super-regularity partition $(R,\mathcal{V})$}	if $G[V_i,V_j]$ is $(\eps,d)$-super-regular for all $ij\in E(R)$.
We sometimes simply say that $G$ admits the regularity partition $(R,\mathcal{V})$ if $\eps$ and $d$ are clear. If $G$ admits a regularity partition $(R,\mathcal{V})$, then $R$ is called \emph{the reduced graph} for $G$.
For given $(R,\mathcal{V})$, let $G[R,\mathcal{V}]= \bigcup_{ij\in E(R)} G[V_i,V_j]$ be the graph we obtain from $G$ after deleting all the edges from non-regular pairs, sparse pairs and the edges within each $V_i$.
\end{defn}
This concept of $\eps$-regularity partition is useful as every graph $G$ admits such a partition $(R,V_1,\dots, V_r)$ with some constant $r$ depending only on $\eps$. This is captured by the celebrated regularity lemma. 
For our purpose, we need the following version with some additional properties.

\begin{lemma}[Szemer\'edi's regularity lemma] \label{regularity}
Let $M,t \in \N$ with ${0 <  1/n\ll \frac{1}{M} \ll \eps, 1/t, 1/s, 1/k \leq 1}$ and $d > 0$. 
Suppose that for each $j\in [2]$, $G^j$ is an $n$-vertex graph on the vertex set $V$ and $\mathcal{U}=(U_1,\dots, U_s)$ is an equitable partition of $V$.
Then there exists a number $r$ which is divisible by $k$ and two graphs $R^1,R^2$ on the vertex set $[r]$, and an equitable partition of $V$ into $\mathcal{V}=(V_1, \ldots, V_r)$ satisfying the following for each $j\in [2]$.
\stepcounter{propcounter}
\begin{enumerate}[label = {\bfseries \emph{\Alph{propcounter}\arabic{enumi}}}]
    \item $t \le r \le M$.
    \item $|V_i| = \frac{n}{r} \pm 1$ for each $i\in [r]$.
    \item $\mathcal{V}$ is a refinement of $\mathcal{U}$. In other words, every set $V_i$ is a subset of $U_j$ for some $j\in [s]$.
    \item For each $j\in [2]$, $G^j$ admits an $(\eps,d)$-regularity partition $(R^j,\mathcal{V})$.
    \item For $G'^j=G^j[R^j,\mathcal{V}]$, all vertices $v\in V(G)$ satisfies $d_{G'^j}(v)\geq d_{G^j}(v) - (d+2\eps)n$.
    \item For some $\delta \ge 0$, if $\delta(G^j)\geq \delta n$, then $\delta(R^j)\geq (\delta -d - 2\eps)r$.\label{cond: reg min deg}
\end{enumerate}
\end{lemma}
This lemma can be proved by following the outline of the proof by Kolm\'{o}s and Simonovits~\cite[Theorem~1.18]{Komlos1996} while they only mentioned an edge-colored graph.

We will often modify the vertex partition while not ruining the $\eps$-regularity too much. The following two lemmas will be useful for this.

\begin{lemma}\label{lem:tiny_part_good_enough}
Suppose $0< \e< \a \le 1$. If $G$ is an $(\eps,d)$-regular graph with the vertex partition $(A,B)$ and $|A'|\ge \a |A|$ and $|B'|\ge \a |B|$, then $G[A',B']$ is $(\frac{\e}{\a},d)$-regular.
\end{lemma}

\begin{lemma}\label{lem:adding_tiny_parts}
Suppose $0< \eps \ll d \le 1$. Suppose that $(A,B)$ is an $(\eps,d)$-regular pair in $G$ and $A,B,A',B'$ are pairwise disjoint sets of vertices in $G$.
If $|A'|\le \e |A|$, $|B'|\le \e |B|$, then $(A\cup A',B\cup B')$ is a $(4\e^{1/2},d)$-regular pair in $G$.
\end{lemma}

\begin{lemma}\label{lem:union_regular}
Let $0<\frac{1}{n} \ll \e \ll d$.
Let $A, B_1,\dots, B_\ell$ be pairwise disjoint sets in a graph $G$ with $|A|, |B_i| \geq n$.
If $(A,B_i)$ is $(\eps,d+)$-regular pair in $G$, then $G$ contains a subgraph $G'$ such that $(A,\bigcup_{i\in [\ell]} B_i)$ is an $(\eps^{1/3},d+)$-regular pair in $G'$.
\end{lemma}

The following lemma converts an $(\eps,d+)$-regular pair into an $(\eps,d)$-super-regular pair.

\begin{lemma}\label{lem:large_super_regular}
Suppose $0<  \e\ll d< 1$. If $(A,B)$ is $(\e,d)$-regular pair in $G$, then there exists $A'\subseteq B, B'\subseteq B$ such that $|A'|\ge (1-\e)|A|$, $|B'|\ge (1-\e)|B|$, and $(A',B')$ is $(2\e, d)$-super-regular pair in $G$.
\end{lemma}

\begin{lemma}\label{lem: nice_regularity}
Suppose $0<1/n \ll 1/r \ll \eps_0 \ll \eps, d, 1/k < 1$.
    Suppose that $\delta(G)\geq (1- \frac{1}{k}+\eps)n$ and $G$ admits an $(\eps_0,d)$-regularity partition $(R,\mathcal{V})$ with $\mathcal{V}=(V_1,\dots, V_r)$ such that $|V_i|= (1\pm \eps_0)\frac{n}{r}$ and $k$ divides $r$. Suppose that $R$ contains a $K_k$-factor $Q$. Then, there exists a new partition $\mathcal{V}'=(V'_1,\dots,V'_r)$ of $V(G)$ satisfying the following.
    \begin{enumerate}
    \item[(a)] $|V'_i| = \left(1\pm \eps_0^{1/2} \right)n/r$ for all $i \in [r]$,
    \item[(b)] for each $ij \in E(Q)$, the graph $G[V'_i,V'_j]$ is $(\eps_0^{1/2},d)$-super-regular, and
    \item[(c)] for each $ij \in E(Q)$, we have $||V'_i|-|V'_j||\le 1$.
    \item[(d)] $G$ admits an $(\eps_0^{1/2},d)$-regularity partition $(R,\mathcal{V}')$.
    \end{enumerate}
\end{lemma}
\begin{proof}
Let $Q^1,\dots, Q^{r/k}$ be the cliques in $Q$.
    For each $ij\in E(Q)$, apply Lemma~\ref{lem:large_super_regular} to the pair $(V_i,V_j)$, and delete more vertices if necessary, to obtain subsets $V^*_1,\dots, V^*_r$ of $V_1,\dots, V_r$, respectively, such that the following holds:
\begin{enumerate}
    \item $|V^*_i|=(1-2k\eps_0)\frac{n}{r}$,
    \item for all $ij\in E(Q)$, the pair $(V^*_i, V^*_j)$ is $(k\eps_0,d)$-super-regular.
\end{enumerate}
This yields $|V\setminus V^*_{[r]}| \leq k\eps_0 n.$
For each vertex in $V\setminus V^*_{[r]}$, we know that $d_G(v)\ge \left(1-\frac{1}{k}+\eps\right)n$ holds, so there exists at least $\left(1-\frac{1}{k}+\frac{\eps}{2}\right)r$ indices $i\in [r]$ such that $|N_G(v)\cap V^*_i|\ge \frac{1}{4}\eps |V^*_i| \geq 2d |V^*_i|$. Indeed, if not, then we have 
$$d_G(v) \leq \left(1-\frac{1}{k}+\frac{\eps}{2}\right)r \cdot \left(\frac{n}{r}+1\right) + \frac{1}{4}\eps n < \left(1-\frac{1}{k}+\eps\right)n, $$
a contradiction.
Hence, for each $v \in V\setminus V^*_{[r]}$, there exists at least $\eps r/2$ choices of $j$ such that $v$ has at least $2d|V^*_i|$ neighbors in $V^*_i$ for all $i\in Q^j$. Hence, moving $v$ to one of $V^*_i$ does not ruin the degree condition for super-regularity of $(V^*_i, V^*_{i'})$ with $ii'\in E(Q^j)$.
By checking such possibility for all $v\in V\setminus V^*_{[r]}$ and distributing them in an appropriate manner, it is routine to check that we can distribute the vertices in $V\setminus V^*_{[r]}$ to $V^*_1,\dots, V^*_r$ to obtain the sets $V_1',\dots, V_r'$ satisfying (a) and (c). Then (b) and (d) follow from Lemma~\ref{lem:tiny_part_good_enough}.
\end{proof}

The following lemma allows us to pick many typical pairs with large common neighborhoods.
\begin{lemma}\label{lem: typical pairs}
    Suppose $0<\eps \ll d \leq 1$.
    Let $V_1,\dots, V_k$ be $k$ vertex sets of a graph $G$ such that $(V_i,V_j)$ is $(\eps,d)$-regular in $G$ for every distinct $i,j\in [k]$.
    Then at least $(1-4k^2\eps^{1/2})|V_1||V_2|$ pairs $(u,v)\in V_1\times V_2$ satisfy $|N_G(\{u,v\}) \cap V_j| \geq (d^2-4\eps^{1/2})|V_j|$ for all $j\in [k]\setminus \{1,2\}$ as well as
    $|N_G(u) \cap V_2| \geq (d-4\eps^{1/2})|V_2|$ and  $|N_G(v) \cap V_1| \geq (d-4\eps^{1/2})|V_1|$.
\end{lemma}
\begin{proof}
By Lemma~\ref{lem:large_super_regular}, there are at least $(1- 2k\eps)|V_1|$ choices of $u\in V_1$ having at least $(d-2\eps)|V_j|$ neighbors in $V_j$ for all $j\in [k]\setminus \{1\}$.
For such a choice $u$, consider $U_j = N_{G}(u)\cap V_j$ for each $j\in [k]\setminus \{1\}$ and let $U_1=V_1$. Lemma~\ref{lem:tiny_part_good_enough} ensures that $(V_2,U_j)$ is $(\eps^{1/2},d)$-regular in $G$. Again Lemma~\ref{lem:large_super_regular} ensures that there are at least $(1-2k\eps^{1/2})|V_2|$ choices of $v\in V_2$ such that $|N_G(v)\cap U_j|\geq (d-2\eps^{1/2})|U_j|$ for all $j\in [k]\setminus \{2\}$.
This yields at least $(1-2k\eps)(1-2k\eps^{1/2})|V_1||V_2|\geq (1- 4k^2 \eps^{1/2})|V_1||V_2|$ pairs $(u,v)$ such that for every $j\in [k]\setminus \{1,2\}$, we have $$|N_G(\{u,v\})\cap U_j|\geq (d-2\eps^{1/2})(d-2\eps)|V_j| \geq (d^2-4\eps^{1/2})|V_j|.$$
By the choices of $u$ and $v$, we also obtain 
$|N_G(u) \cap V_2| \geq (d-4\eps^{1/2})|V_2|$ and  $|N_G(v) \cap V_1| \geq (d-4\eps^{1/2})|V_1|$.
\end{proof}

The $\eps$-regularity is a global condition, so it is sometimes difficult to deal with. The following local condition is almost equivalent to $\eps$-regularity, and it is often helpful to us. 

\begin{defn} [Quasi-random]
A bipartite graph with vertex partition $(A,B)$ is $(\eps,p)$-quasi-random if for any two distinct vertices $u,v\in A$, we have $d(u)=(1\pm \e)p|B|$ and $d(u,v)=(1\pm \e)p^2|B|$, and for any two distinct vertices $u,v\in B$, we have $d(u)=(1\pm \e)p|A|$ and $d(u,v)=(1\pm \e)p^2|A|$.
\end{defn}

\begin{theorem}[\cite{Duke1995}]\label{equivalent_to_regularity}
Suppose $1/n \ll \e \ll \a,p \le 1$. Suppose $G$ is a bipartite graph with vertex partition $(A,B)$ such that $|A|=n$, $\a n\le |B|\le \a^{-1}n$, and at least $(1-5\e)n^2/2$ pairs $u,v \in A$ satisfy $d(u),d(v) \ge (p-\e)|B|$ and $d(u,v)\le(p+\e)^2|B|$. Then $G$ is $(\e^{1/6},p)$-regular. In particular, if $G$ is $(\e,p)$-quasi-random, it is $(\e^{1/6},p)$-super-regular.
\end{theorem}

Using this theorem, one can prove the following theorem stating that, for an $\e$-regular pair $(A,B)$, random subsets of $A$ and $B$ somewhat inherit the $\eps$-regularity.

\begin{lemma}\label{lem:randomset_reg}
 Let, $0< 1/n \ll 1/r \ll \varepsilon \ll d\le 1$ and $0<1/n \ll \beta \le 1$.
 Let $G$ be a graph with the equitable $(\eps ,d)$-regular partition $(R,\mathcal{V})$ with $\mathcal{V}=(V_1, \ldots, V_r)$ and $(\eps ,d)$-super-regular partition $(R',\mathcal{V})$ for some $R' \subseteq R$.
 Let $\beta n \leq m \leq (1-\beta)n$ be an integer and $U$ be a random set of vertices of size $m$ chosen uniformly at random and let $U_i=U \cap V_i$, $U'_i = V_i \setminus U$.
 Then with probability at least $1-o(1/n)$, we have 
 \begin{itemize}
    \item $|U_i| =(1 \pm \varepsilon^{1/4}) \frac{m}{r}$ \label{cond:randomset1}
    \item If $ij \in E(R)$, then $(U_i, U_j)$, $(U_i, U'_j)$ and $(U'_i, U'_j)$ are $(\eps^{1/10} ,d+)$-regular pairs \label{cond:randomset2}
    \item If $ij \in E(R')$, then $(U_i, U_j)$, $(U_i, U'_j)$ and $(U'_i, U'_j)$ are $(\eps^{1/10} ,d+)$-super-regular pairs
 \end{itemize}
\end{lemma}
\begin{proof}
 For the first part, by Lemma~\ref{lem:concentration}, the probability that the size of $U_i$ is not $\left(1 \pm \varepsilon^{1/4}\right) \frac{m}{r}$ is at most $$2\exp\left(-(\varepsilon^{1/4}-\varepsilon^{1/2})^2 m^2/8r^2 n \right) \leq o(n^{-2}).$$
The union bound yields that, with probability $1-o(n^{-1})$, the first condition hold for all $i\in [r]$.

 For the second and third parts, suppose that $(V_i,V_j)$ is an $(\eps, d')$-regular pair for some $d'=d'(ij) \geq d$.
By definition, all but $2\eps |V_i|$ vertices in $V_i$ have $(d' \pm \eps)|V_j|$ neighbors in $V_j$ and for each such vertex $v \in V_i$, all but $2{d'}^{-1}\eps |V_i|$ vertices $u \in V_i$ satisfy that $|N(v) \cap N(u) \cap V_j| = (d' \pm 2{d'}^{-1}\eps) (d' \pm \eps) |V_j| = ({d'}^2 \pm 10{d'}^{-1} \eps) |V_j|$.
 For such a fixed vertex $v$, the probability that $|N(v) \cap U_j| = (d' \pm 2\varepsilon)|V_j| \frac{m}{n}$ is at least $$ 1- 2\exp\left(-\varepsilon^2 |V_j|^2 m^2 / 8n^3\right) = 1-o(n^{-4})$$ 
 by Lemma~\ref{lem:concentration}.
 Similarly, for $u, v \in V_i$ with $|N(v) \cap N(u) \cap V_j| = \left({d'}^2 \pm 10{d'}^{-1} \eps\right) |V_j|$, we have 
 $$\mathbb{P}[|N(v) \cap N(u) \cap U_j| = \left({d'}^2 \pm 20{d'}^{-1} \eps\right) |V_j| \frac{m}{n}] \geq 1-2\exp\left(-10^2{d'}^{-2}\varepsilon^2 |V_j|^2 m^2 / 8n^3\right) = 1-o(n^{-4}). $$
 By the union bound, with probability at least $1-o(n^{-3})$, we have 
 $$|N(v) \cap U_j| = \left(d' \pm 2\varepsilon\right)|V_j| \frac{m}{n} \text{ and } |N(v) \cap N(u) \cap U_j| = \left({d'}^2 \pm 20{d'}^{-1} \eps\right) |V_j| \frac{m}{n}$$ 
 holds for at least $(1-2\eps) |V_i|$ vertices $u\in V_i$ and  for at least $(1-2\eps)(1-2{d'}^{-1}\eps)|V_i|^2 \geq (1-4d^{-1} \varepsilon)|V_i|^2$ pairs $u, v \in V_i$.
 Therefore, with probability $1-o(n^{-1})$, the same conclusion holds for every $i, j \in [r]$ with $ij \in E(R)$.
 By Lemma~\ref{equivalent_to_regularity}, we have $(U_i, U_j)$ is an $(\eps^{1/10} ,d+)$-regular pair for every $ij \in E(R)$.
 Furthermore, if $ij \in E(R')$, then $|N(v) \cap U_j| \geq (d'(ij) - 2\varepsilon)|V_j| \frac{m}{n}$ for every $v \in V_i$ by the same argument. Thus, $(U_i, U_j)$ is an $(\eps^{1/10} ,d+)$-super-regular pair for every $ij \in E(R')$.
 
 As $V \setminus U$ is a uniformly chosen random set of size $n-m \in [\beta n, (1-\beta)n]$, the same argument shows that $(U_i, U'_j)$ is an $(\eps^{1/10} ,d+)$-regular pair for every $ij \in E(R)$ and $(U_i, U'_j)$ is an $(\eps^{1/10} ,d+)$-super-regular pair for every $ij \in E(R')$. The same holds for the pairs $(U'_i, U'_j)$.
 This finishes the proof of the lemma.
\end{proof}

\subsection{The blow-up lemma}

The following blow-up lemma provides another reason why $\eps$-regularity is important. It essentially says that if a bounded degree graph $H$ admits a vertex partition $(R,X_1,\dots, X_r)$ and $G$ admits an $\eps$-regularity partition $(R,V_1,\dots, V_r)$, then we can embed $H$ into $G$ as long as each $X_i$ is not larger than~$V_i$.

\begin{lemma}[The blow-up lemma] \label{lem:technical_blow} 
Suppose $0 < 1/n \ll \g \ll \eps \ll \beta \ll \a,d,1/\D,1/k,1/s \le 1$ and $0< 1/n \ll 1/r\leq 1$. Suppose that $G,H$ are graphs and $\M$ is a multi-hypergraph on the vertex set $V(H)$, $R$ is a graph on vertex set $[r+s]$, $\mathcal{V}=(V_1,\dots, V_{r+s})$ is a partition of $V(G)$, $\mathcal{X}=(X_1,\dots, X_k)$ is a partition of $V(H)$ and $\alpha n/r \leq |V_i|\leq \alpha^{-1} n/r$ for all $i\in [k]$.
Let $g:[k]\rightarrow [r]$ be a function and $f:E(\M)\rightarrow [r+1,r+s]$ be a function and $T\subseteq V(H)$ be a set of at most $\g n/r$ vertices. For all $i\in [k]$ and $x\in T\cap X_i$, we have a vertex set $A_x\subseteq V_i$. Suppose that the following hold.
\stepcounter{propcounter}
\begin{enumerate}[label = {\bfseries \emph{\Alph{propcounter}\arabic{enumi}}}]
    \item 
The graph $G$ admits $(\eps,d)$-super-regularity partition $(R,\mathcal{V})$ and $H$ admits a vertex partition $(R',\mathcal{X})$ with $V(R')=[k]$.
Also, $g:V(R')\rightarrow V(R)$ is a graph homomorphism from $R'$ to $R$, i.e. $g(i)g(j)\in E(R)$ for all $ij\in E(R')$. \label{blowup_assumption_0}
\item For each $i\in [r]$, we have $\sum_{j\in g^{-1}(i)} |X_j| \leq |V_i|$. Moreover, 
if $\sum_{j\in g^{-1}(i)} |X_j| \leq (1-\eps^{1/2})|V_i|$ for all $i$, then $G$ may only admit $(\eps,d)$-regularity partition $(R,\mathcal{V})$ instead of admitting super-regularity partition.
\label{blowup_assumption_1}
    \item $\D(H) \le \D$, $\D(\M) \le \D$, and the rank of $\M$ is at most $\D$, \label{blowup_assumption_2}
    \item For all $i\in [r+1,r+s]$, $j \in [k]$, and $e\in E(\M)$, if $f(e)=i$ and $e\cap X_j \neq \varnothing$, then $ij\in E(R)$.  \label{blowup_assumption_3}
    \item For all $i\in [k]$ and $x\in T\cap X_i$, we have $|A_x|\ge \a|V_i|$. Moreover, the vertices in $T$ does not belong to any edges in $\M$.\label{blowup_assumption_4}
\end{enumerate}
Then there exists an embedding $\phi$ of $H$ into $G$ such that $\phi(X_i) \subseteq V_{g(i)}$ for all $i\in [k]$ and
\stepcounter{propcounter}
\begin{enumerate}[label = {\bfseries \emph{\Alph{propcounter}\arabic{enumi}}}]
    \item for all $i\in [k]$ and $x\in X_i$, we have $\phi(x)\in V_i$, and for each $x\in T$, we have $\phi(x)\in A_x$, \label{blowup_conclusion_1}
    \item for each $e\in E(\M)$, we have $|N_G(\phi(e))\cap V_{f(e)}| \ge \beta |V_{f(e)}|$. \label{blowup_conclusion_2}
\end{enumerate}
\end{lemma}
Note that the usual version of the blow-up lemma is stated with the sets $X_1,\dots, X_k$ and $V_1,\dots, V_k$ and assuming that $g:[k]\rightarrow [k]$ is the identity function, or more generally, when $g$ is an injective function. However, the above statement is actually equivalent with it.
Indeed, by taking an appropriate random partition of $V_i$ into $\{ V'_j: j\in g^{-1}(i)\}$ with $|V'_j|\geq |X_j|$ and applying Lemma~\ref{lem:randomset_reg}, it is straightforward to see that the above version with the function $g$ is equivalent with the original blow-up lemma. Introduction of the function $g$ will be convenient for later applications.

There are two more additional features in the above lemma compared to the usual version of the blow-up lemma: (1) There is a small number of vertices $x$ which we have to embed into a set $A_x$.
(2) There is a hypergraph $\M$ on a small vertex set $V(\M)\subseteq V(H)$ where the images of an edge $e\in \M$ must be embedded so that they have a large common neighborhood in $V_{f(e)}$. 


Imagine $V(H)=X_0\cup X_1\cup X_2$ and $H_0=H[X_0]$ is already embedded by $\phi_0$ into $G$ and we wish to embed $H'= H[X_1]$ now while there are remaining part $H''=H[X_2]$ which we will embed later.
As $H_0$ is already embedded, if some vertices $x\in V(H')$ has neighbors $y,z$ in $X_0$ we have to embed $x$ into the common neighborhood $N_G(\phi_0(y),\phi_0(z))$. This information is encoded by the target set $A_x$.
Since we want to later embed a vertex $x\in X_2$, if $x$ has neighbors $y,z$ in $X_1$, we want to embed $y,z$ now so that the common neighborhood $N_G(\phi(y),\phi(z))\cap X_2$ is large. This information about common neighborhoods of the vertices in $X_2$ is encoded by the hypergraph $\M$ and the function $f$, and \ref{blowup_conclusion_2} ensures the common neighborhood $N_G(\phi(y),\phi(z))\cap X_2$ to be large so that we can apply this lemma again later.

Since these two features are standard by now (for example, Lemma~5.1 in \cite{Kim2019}
yields these two features with an additional assumption that all $V_i$ have almost the same sizes; however, this extra assumption is not necessary to obtain these two features), we omit the proof by just noting that these two additional features can be verified by almost exactly following the proof of the blow-up lemma in~\cite{Rodl1999}.

\subsection{Graph systems}
Assume that we have a collection $\G=(G_1,\dots, G_h)$ of graphs with the same vertex set $V$. 
We consider each index in $[h]$ as a color, and $G_i$ as the graph with edges colored with $i$. For each $uv\in \binom{V}{2}$, let $\Lambda_\G(uv)=\{ i\in [h]: uv\in E(G_i)\}$ be the set of colors $i$ where $uv$ belong to the graph $G_i$. 
For a set $C\subseteq [h]$ of colors, let $\G_C=( G_i: i\in C)$. 
For $\eta\in [0,1]$, by only considering edges belonging to at least $\eta |C|$  graphs in $\G_C$, we obtain the following graph which will be useful to us. 
$$\G_C^\eta = \{ uv: |\Lambda_\G(uv)\cap C|\geq \eta|C| \}.$$

We call $\G^{\eta}=\G^{\eta}_{[h]}$ the \emph{$\eta$-fraction graph} of the graph system $\G$.
This concept of fraction graph was introduced in \cite{Montgomery2022} by Montgomery, M\"{u}yesser and Pehova. They observed that this concept is very useful for finding graph partial transversals.
Indeed, once we find an embedding $\phi$ of $H'\subseteq H$ into $\G^\eta$ and $H'$ has at most $\eta h$ edges, then we can greedily choose a color for each edge $uv\in E(\phi(H'))$ from $\Lambda_{\G}(uv)$ in such a way that all edges have different colors. In addition, the following lemma shows that the fraction graph somehow inherits the minimum degree condition.

\begin{lemma}[{\cite[Proposition~2.1]{Montgomery2022}}] \label{lem_mindeg}
Suppose $0<\eta, \delta \leq 1$.
Let $\G=(G_1,\dots, G_h)$ be a collection of $n$-vertex graphs with $\delta(G_i)\geq \delta n$ for all $i\in [h]$.
Then we have $\delta(\G^{\eta}) \geq (\delta - \eta) n$.
\end{lemma}

Hence, an $\eta$-fraction graph with small $\eta$ is useful for graph embedding. However, this only lets us embed small bits of $H$. Even if we repeat this process of embedding, a significant number of vertices of $H$ will remain unembedded. The following color-absorption lemma allows us to be able to choose colors for some edges all at once at the end.

\begin{lemma} [\cite{Montgomery2022}] \label{absorption}
Let $0<1/n \ll \g \ll \b \ll \eta < 1$. 
Let $\G=(G_1,\dots, G_h)$ be a collection of graphs with $h\geq \eta n$ and $F\subseteq \G^{\eta}$ be a subgraph of $\G^{\eta}$ with $e(F)=\beta n$.
Then there exist disjoint sets $A, C \subset [h]$, with $|A| = e(F) - \g n$ and $|C| = \b n$ such that the following property holds. Given any subset $C' \subset C$ with $|C'| = \g n$, there is a transversal in $\G_{A\cup C'}$ isomorphic to $F$.
\end{lemma}

Note that this lemma is not stated exactly as above in \cite{Montgomery2022}. For completeness, we deduce it from the lemma below stated in \cite{Montgomery2022}.

\begin{lemma}[{\cite[Lemma~3.3]{Montgomery2022}}]\label{lem:absorber_prelim}
    Let $\alpha \in (0, 1)$, let $n, m$ and $\ell\geq 1$ be integers satisfying $\ell\leq \alpha^7 m/ 10^5$ and $\alpha^2 n \geq 8m$. Let $H$ be a bipartite graph on vertex classes $A$ and $B$ such that $|A| = m$, $|B| = n$ and, for each $v\in A$, $d_H(v)\geq \alpha n$.
    
    Then, there are disjoint subsets $B_0, B_1\subset B$ with $|B_0| = m-\ell$ and $|B_1|\geq \alpha^7 n/ 10^5$, and the following property. Given any set $U\subset B_1$ of size $\ell$, there is a perfect matching between $A$ and $B_0\cup U$ in $H$.
\end{lemma}

\begin{proof}[Proof of Lemma~\ref{absorption}]
Let $K$ be the bipartite graph with vertex classes $E(F)$ and $[h]$, where $(e,i)$ is an edge if $e\in G_i$. Then Lemma~\ref{lem:absorber_prelim} applies with $(\ell,m,n,\alpha)=(\gamma n,\beta n, h, \eta )$. This yields disjoint $A,C\subseteq[m]$ with $|A|=e(H)-\gamma n$ and $|C|\ge 10\beta n$, such that, for any $C'\subseteq C$ of size $\gamma n$ there is a perfect matching between $E(F)$ and $A\cup C'$. Such a matching corresponds to a $\G_{A\cup C'}$-transversal isomorphic to $F$, as required.
By taking a subset of $C$ with size exactly $\beta n$ to be our set $C$, we prove the lemma.
\end{proof}

\section{Proof sketch}\label{sec: proof sketch}

For the proof of Theorem~\ref{main}, we develop the ideas from \cite{Montgomery2022}, where Montgomery, M\"uyesser and Pehova laid out a general strategy to find a $\mathcal{G}$-transversal isomorphic to a spanning graph $H$ in the collection $\mathcal{G}=(G_i: i\in [h])$. 
Using low bandwidth of $H$, we first partition $H$ into vertex-disjoint induced subgraphs $H_0,\dots, H_{m+1}$ so that $e(H_0)$ is much bigger than $e(H_i)$ for any $i\in [m+1]$ and
most of the edges of $H$ are within $H_i$ and the remaining very small amount of edges are all between two consecutive vertex sets $V(H_{i-1})$ and $V(H_{i})$.
\newline

\noindent {\bf Step 1. Find a color absorber {\boldmath $(\phi(H_0),A,B,C)$}.} We find an embedding $\phi$ of $H_0$ into $\mathcal{G}$ while finding $A,C\subseteq [h]$ such that the following holds: for any given $C'\subseteq C$ with $|C'|=e(H_0)-|A|$, the graph $\phi(H_0)$ can be colored using exactly the colors in $A\cup C'$. Let $B= [h]\setminus(A\cup C)$. \newline

\noindent {\bf Step 2. Embed {\boldmath $H_1,\dots, H_{m+1}$} using all the colors in {\boldmath $B$} and some colors in {\boldmath $C$} while not using any colors in {\boldmath $A$}.}
\begin{itemize}
    \item[]{\bf Step 2-1.} Extend the embedding of $H_0$ into the embedding of $H_0\cup \dots \cup H_{m}$, and color all edges not in $E(H_0)$ using colors in $B\cup C$ while making sure that most of the colors in $B$ are used.
    \item[]{\bf Step 2-2.} Extend the embedding by embedding the last piece $H_{m+1}$ and color the edges so that all colors in $B$ are used, and no colors in $A$ are used. We make sure that all vertices are covered as well.
\end{itemize}

\noindent {\bf Step 3. Color the edges in {\boldmath $H_0$} exhausting all the remaining colors.} Take the remaining unused colors in $C'$ and apply the color absorbing property of $H_0$ to color the edges of $H_0$ with the colors in $A\cup C'$, finishing a rainbow coloring of $H$. \newline

One of two innovative ideas of \cite{Montgomery2022} that makes this scheme work is the use of color absorber.
This allows one to put aside a small part $H_0$ of the graph which can be used to finish the rainbow coloring at the end while giving an extra room to embed the remaining graph. 
The other one is the fraction graph. We break the graph $H$ into small pieces $H_0,\dots, H_{m+1}$.
As we do not color $H_0$ until the end, at the time of embedding $H_i$, we have at least $e(H_0)-|A|+ e(H_i)$ colors available which is much more than the number $e(H_i)$ of colors we need. 
Hence, by taking an appropriate fraction graph and using Lemma~\ref{lem_mindeg}, we can reduce our problem to finding a copy of $H_i$ within one graph $\mathcal{G}^\eta_{B'}$ for some $\eta>0$ and some set $B'$ of colors.

In order for this scheme to work, we need to achieve two goals. 
First, we need to be able to embed each piece $H_i$ while making sure that some set $B$ of colors are all exhausted, and all vertices are covered. This is necessary to carry out Step 2-2. 
Second, we need to be able to establish the connections between two pieces $H_{i-1}$ and $H_i$. 
Building connection was relatively simpler in \cite{Montgomery2022}, where they embedded an $F$-factor and trees. 
For an $F$-factor, such a connection is unnecessary since an $F$-factor can be partitioned into small pieces without any edges between two distinct pieces. 
For trees, using $1$-degeneracy, we only need to keep track of how one vertex is embedded in the previous embedding, and establish the connection from there. However, for a general graph $H$, providing such a connection turns out to be much more complicated.

To achieve our first goal, we need a tool for embedding low-bandwidth graph $H$ into the fraction graph $G= \mathcal{G}^{\eta}_{B'}$ with large minimum degree for some color subset $B'\subseteq [h]$. Lemma~\ref{lem:technical_bandwidth} allows us to embed each part $H_i$ in such a way that the following holds.
\begin{enumerate}
    \item[$\bullet$] A small number of vertices $x$ are embedded into a pre-designated set $A_x$.
    \item[$\bullet$] Specific sets (edges of a hypergraph $\M$) of the vertices in $H_i$ are embedded into sets of vertices whose common $G$-neighborhoods are large. 
\end{enumerate}
In the applications, the set $A_x$ will be a common $G$-neighborhood of pre-embedded images of the $H$-neighbors of the vertex $x\in V(H_{i})$ and the set $e\in \M$ will be a set of $H$-neighbors of a vertex $y\in V(H_{i+1})$ in $V(H_i)$. The second bullet point for $H_i$ will provide the set $A_y$ for the next round of embedding for $H_{i+1}$. These two properties will ensure the continuation of our embedding scheme. Moreover, the first bullet point later will be useful in Step~2-2 to ensure that we use up all the colors in $B$. This can be done by hand-picking a $3$-independent matching $M$ in $H_{m+1}$ and a rainbow matching in $\mathcal{G}$ using all the remaining colors in $B$. Once embedding the matching $M$ in $H_{m+1}$ into the rainbow matching in $\mathcal{G}$ and setting up appropriate sets $A_x$ for each $H$-neighbors of the vertices in $M$, application of Lemma~\ref{lem:technical_bandwidth} will yield a desired embedding so that we can use up all the colors in $B$.

At first glance, it seems that this already provides enough connections between $H_i$ and $H_{i+1}$ for our purpose. However, one key difficulty is that the graph we use to embed $H_i$ is different from the graph we use to embed $H_{i+1}$. Namely, once we embed $H_i$ into $\mathcal{G}^{\eta}_{C^i}$ for some color set $C^i$ and color each edge, we have to discard the used colors. Then we have to embed $H_{i+1}$ into a new fraction graph $\mathcal{G}^{\eta}_{C^{i+1}}$ where $C^{i+1}$ is a set different from $C^{i}$. This might not be complicated if $C^{i+1}$ is slightly different from $C^{i}$, but we need to deal with the cases where two sets are substantially different.

Our plan to overcome this difficulty is roughly as follows. We use Lemma~\ref{lem:technical_bandwidth} to embed $H_i$ into $G^i= \mathcal{G}^{\eta}_{C^i}[U^i]$ for some specific subset $U^i$ of $V(\mathcal{G})$ and specific color set $C^{i}\subseteq[h]$ using the regularity structure $(R^i, \mathcal{V}^i)$ of $G^i$.
Similarly, we wish to embed $H_{i+1}$ into $G^{i+1}=\mathcal{G}^{\eta}_{C^{i+1}}[U^{i+1}]$ using a new regularity structure $(R^{i+1}, \mathcal{V}^{i+1})$ of $G^{i+1}$.
In order to build connections between two distinct regularity structures,  we apply the regularity lemma (Lemma~\ref{regularity}) with additional properties to find another common regularity partition $(R,\mathcal{V})$ of both graphs $G^i$ and $G^{i+1}$ which almost refines both $(R^i, \mathcal{V}^i)$ and $(R^{i+1}, \mathcal{V}^{i+1})$. 
By finding an appropriate clique-walk within this finer regularity structure $R$ and exploiting some technical properties, we can establish a desired connection between two regularity partitions. This connection allows us to smoothly transit the embedding of $H_i$ using $G^i$ into the embedding of $H_{i+1}$ using $G^{i+1}$. The connection lemma~(Lemma~\ref{lem: connection}) provides such a smooth transition.

\section{Connection Lemma}
As it was explained in the previous section, we wish to establish smooth connections between two graphs $G^1, G^2$ each having a regularity partition over a common vertex set. 
In other words, we need to take some initial segments of the bandwidth ordering of $H_{i+1}$ (as introduced in the last section) and embed it into $G^1\cup G^2$ in such a way that it provides the connection between $H_i$ and $H_{i+1}$. The graph $H$ in the following lemma will be the subgraph of $H_{i+1}$ induced by the initial segments of the ordering, which has much smaller number of vertices than $G^1$ or $G^2$.

\begin{lemma}[Connection Lemma]\label{lem: connection}
    Let $0 \ll \frac{1}{n} \ll \alpha \ll \frac{1}{r}, \eps_0 \ll  d', \beta \ll d  \ll \eps \ll \frac{1}{k}, \frac{1}{\Delta} \le 1$ and $0\ll \alpha \ll \frac{1}{\ell} \le 1$.
    Let $H$ be a graph with $\Delta(H)\leq \Delta, \chi(G)=k$ and an $\alpha n$-bandwidth partition $W_1,\dots, W_{\ell}$ with $\ell \geq 4k^3$ and $k$ divides $\ell$. 
   Let $G^1, G^2$ be two graphs on a common vertex set $V$ of size $n$.
   For each $j\in [2]$, we have  $\delta(G^j) \geq (1-\frac{1}{k}+\eps)n$ and $G^j$ admits an $(\eps_0,d)$-regularity partition $(R^j, \mathcal{V}^j)$ with $\mathcal{V}^{j}= (V^j_1,\dots, V^j_{r_j})$ and $r_1,r_2\leq r$ and $|V^j_i| = (1 \pm \eps^{1/2} )\frac{n}{r_j}$ for every $i\in [r_j]$. Let $[k]$ be a clique in both $R^1$ and $R^2$. 
   
   For each $i\in [k]$ and $x\in W_{i}$, we have a target set $A_x\subseteq V^1_i$ for each $x$ with $|A_x|\geq d'|V^1_i|$. Let $\M$ be a multi-hypergraph on the vertex set $W_{[\ell -k +1,\ell]}$ with rank at most $\Delta$ and $\Delta(\M)\leq \Delta$ and a function $f:E(\M)\rightarrow [k-1]$ is given such that every hyperedge $e\in \M$ satisfies $e\subseteq W_{[\ell-k+f(e)+1,\ell]}$.
    Then there exists an embedding $\varphi$ of $H[W_{[\ell]}]$ into $G^1\cup G^2$ such that the following hold.
\stepcounter{propcounter}
\begin{enumerate}[label = {\bfseries \emph{\Alph{propcounter}\arabic{enumi}}}]
        \item For each $x\in W_{[k]}$, we have $\varphi(x) \in A_x$. \label{eq: connection conclusion 1}
        \item For each $j \in [k-1]$ and $e\in \M$ with $f(e)=j$, we have 
        $$\left|\bigcap_{x\in e} N_{G_2}(\varphi(x)) \cap (V_{j}^2 \setminus \varphi(W_{[\ell]}))\right| \geq \beta \left|V_{j}^2 \setminus \varphi(W_{[\ell]})\right|.$$ \label{eq: connection conclusion 2}
    \end{enumerate}
\end{lemma}
\begin{proof}
We choose further constants $\eps_1, \hat{d}, \hat{r}_0, \hat{r}_1$ so that we have 
$$ 0< \frac{1}{n} \ll \alpha \ll \frac{1}{\hat{r}_1} \ll \frac{1}{\hat{r}_0} \ll \eps_1 \ll \frac{1}{r}, \eps_0 \ll d', \beta \ll \hat{d} \ll d < 1.$$
In order to embed $H$ into $G^1\cup G^2$, we want to find a finer regularity partition that both $G^1$ and $G^2$ admit. For this, take the common refinement $\mathcal{U} = \{ V^1_i\cap V^2_j : i\in [r_1], j\in [r_2]\}$ of $\mathcal{V}^1$ and $\mathcal{V}^2$.
For each set $U$ in $\mathcal{U}$, we partition $U$ into sets of size exactly $\frac{n}{r^5}$ and one set of less than $\frac{n}{r^5}$ remaining vertices. We collect all the remaining vertices to obtain a set $U_0$ with $|U_0|\leq \frac{n}{r^5}\cdot r^2 \leq \frac{n}{r^3}$. Collect the obtained sets of size $\frac{n}{r^5}$ to get a partition $\mathcal{U}'$ of $V'= V\setminus U_0$ into sets of size $\frac{n}{r^5}$.

We apply Lemma~\ref{regularity} to the graph $G^1[V'], G^2[V']$ and the regularity parameter $\eps_1/4$ and $\hat{d},\hat{r}_0, \hat{r}_1$ playing the roles of $d,t, M$, respectively, we obtain an equitable partition $\widehat{\mathcal{V}} = (\widehat{V}_1,\dots \widehat{V}_{\hat{r}})$ satisfying the following.
\stepcounter{propcounter}
\begin{enumerate}[label = {\bfseries \emph{\Alph{propcounter}\arabic{enumi}}}]
        \item $\hat{r}_0 \leq \hat{r}\leq \hat{r}_1$, hence $\alpha \ll \frac{1}{\hat{r}} \ll \eps_1$.
        \item $|\widehat{V}_i| = \frac{1}{\hat{r}}|V'| \pm 1$ for each $i\in [\hat{r}]$.
        \item For each $j\in [2]$, the set $[\hat{r}]$ partitions into $\mathcal{I}^j=(\hat{I}^j_1,\dots, \hat{I}^j_{r_j})$ such that for each $i\in [r_j]$, we have $\widehat{V}_{I^j_i} \subseteq V^j_i$ and $|\widehat{V}_{I^j_i}|\geq (1-\frac{2}{r^3})|V^j_i|$. \label{eq: common reg 3}
        \item There exist two graphs $\hat{R}^1$ and $\hat{R}^2$ such that for each $j\in [r]$, the graph $G^j[V']$ admits $\left(\eps_1,\hat{d}\right)$-regularity partition $(\hat{R}^j, \widehat{\mathcal{V}})$.\label{eq: common reg 4}
        \item For each $j\in [2]$ and $i\in [\hat{r}]$, the number of $i'\in [\hat{r}]$ such that $(\widehat{V}_i,\widehat{V}_{i'})$ is not $\left(\eps_1,\hat{d}\right)$-regular in $G^j$ is at most $\eps_1 \hat{r}$. 
        \label{eq: common reg 5}
        \item For each $j\in [2]$, $\delta(\hat{R}^j)\geq (1- \frac{1}{k}+\frac{\eps}{2})\hat{r}$.
    \end{enumerate}
Five of the above properties are obvious from the conclusion of the Lemma~\ref{regularity}. The inequality in \ref{eq: common reg 3} follows by noting that $|U_0\cap V_i^j|\leq \frac{n}{r^5}\cdot r \leq \frac{2}{r^3} |V^j_i|$ for all $j\in [2]$ and $i\in [r_j]$. 
For each $j\in[2]$, let $s_j = \frac{\hat{r}}{r_j}$, then 
\ref{eq: common reg 3} implies that $|\hat{I}^j_i| = (1\pm 2\eps^{1/2})s_j$ for all $j\in [2]$ and $i\in [r_j]$.

We wish to embed $H$ into $G^1\cup G^2$ using the common regularity partition $\mathcal{V}$. Roughly, we will find a clique-walk $W$ from a clique $K_1$ in $\hat{R}^1$ to a clique $K_2$ in $\hat{R}^2$. In this clique-walk $W$, the edges between the last $k$ vertices belong to $\hat{R}^2$, and the rest belongs to $\hat{R}^1$. Then using the blow-up lemma, we can embed $H$ into this clique-walk, establishing the connection. However, a naive approach like this does not yield the two conditions \ref{eq: connection conclusion 1} and \ref{eq: connection conclusion 2}. 

For \ref{eq: connection conclusion 1} to hold, we need to choose a clique $K^1 = \{i_1,\dots, i_k\}$ so that $\widehat{V}_{i_j}\subseteq V^1_{j}$ and $\widehat{V}_{i_j}$ contains many vertices from $A_x$ for $j\in [k]$ and the vertices $x\in W_j$. 
Moreover, to ensure \ref{eq: connection conclusion 2}, we need to choose a clique $K^2 = \{i'_1,\dots, i'_k\}$ so that $\widehat{V}_{i'_j}\subseteq V^2_{j}$ and 
analyze how many common neighbors the vertices from $\widehat{V}_{i'_j}$ have in $V^2_{j'}$.
Each set $\widehat{V}_{i'_j}$ has at most $\frac{n}{\hat{r}} \pm 1 \ll \eps_0\frac{n}{r}$ vertices. This means that $\widehat{V}_{i'_j}$ is too small to inherit the regularity property of $V^2_j$ even if $\widehat{V}_{i'_j}\subseteq V^2_j$. Thus, all those vertices in $\widehat{V}_i$ could be exceptional vertices in $V^2_j$ with respect to  the original regularity pairs in $\mathcal{V}^2$ (for example, $G^2[\widehat{V}_{i'_j}, V^2_{j'}]$ could even be empty for all $j'\in [k]$), then it is not possible to obtain~\ref{eq: connection conclusion 2}.

We overcome these difficulties by creating a clique-walk transiting from $R^1$ to $\hat{R}^1$, then from $\hat{R}^1$ to $\hat{R}^2$, and finally from $\hat{R}^2$ to $R^2$. In this clique-walk, we want that the initial several sets are from $\mathcal{V}^1$, and next several sets are not-too-small subsets of the sets in $\mathcal{V}^1$, and the next sets are from $\widehat{\mathcal{V}}$ and then subsets of the sets in $\mathcal{V}^2$ and final sets are from $\mathcal{V}^2$.
For this, we need to establish $\eps_0$-regularity between a set from $\mathcal{V}^j$ and a set from $\widehat{\mathcal{V}}$. This can be done by finding a copy of complete $k$-partite graph $K_{1,\dots, 1, m,\dots,m}$ in $\hat{R}^j$ with large $m$. Such each vertex in $K_{1,\dots,1,m,\dots, m}$ corresponds to a vertex set in $\widehat{\mathcal{V}}$, thus by ensuring that the $m$ vertices in the last part yield a large subset of $V^j_i\in \mathcal{V}^j$, we can establish the $\eps_0$-regularity between small sets and big sets as desired.
In order to build desired copies of $K_{1,\dots,1,m,\dots,m}$, we count the number of $K_k$'s within $R^1, R^2, \hat{R}^1, \hat{R}^2$. Both $\hat{R}^j$ and $R^j$ inherit the minimum degree of $G^j$, so they also inherit the density of $K_k$'s within $G^j$. By such analysis, we can find some $k$-clique $K^j$ and some initial segments of the clique-walk in $\hat{R}^1$, which provides good co-degree conditions in $R^j$ and desired the $\eps_0$-regularity.

\begin{claim}\label{cl: clique counting}
Let $j\in [2]$ and  $Q=\{q_1,\dots, q_k\}$ be a copy of $K_k$ in $R^j$.
Then, the $k$-partite subgraph $\hat{R}^j[\hat{I}_{q_1}^j, \hat{I}_{q_2}^j, \ldots, \hat{I}_{q_k}^j]$ has at least $2k\hat{d} s_j^k$ copies of $K_k$.
\end{claim}

\begin{poc}
Assume that $\hat{R}^j[\hat{I}_{q_1}^j, \hat{I}_{q_2}^j, \ldots, \hat{I}_{q_k}^j]$ contains $m$ copies of $K_k$.
As $Q$ forms a copy of $K_k$ in $R^j$, the standard counting lemma yields that the $k$-partite subgraph $G^j[V^j_{q_1}\setminus U_0,\dots, V^j_{q_k}\setminus U_0]$ of $G^j$ contains at least 
$$ (d-\eps_0)^{\binom{k}{2}} \prod_{i\in [k]} \left|V^j_{q_i}\right| - |U_0| \sum_{l\in [k]} \prod_{i\in [k]\setminus \{l\}} \left|V^j_{q_i}\right| \geq \left(\frac{d}{2}\right)^{\binom{k}{2}} \left(\frac{n}{r_j}\right)^k$$ 
copies of $K_k$.

We now count these copies of $K_k$ in $G^j[V^j_{q_1}\setminus U_0,\dots, V^j_{q_k}\setminus U_0]$ again in another way. For each copy $K_k = (j_1,\dots, j_k)$ in $\hat{R}^j[\hat{I}_{q_1}^j, \dots, \hat{I}_{q_k}^j]$, the corresponding $k$-partite subgraph $G^j[ \widehat{V}_{j_1},\dots, \widehat{V}_{j_k}]$ contains at most
$\left(\frac{n}{\hat{r}}+1\right)^k $ copies of $K_k$. 
As there are $m$ such copies of $K_k$ in $\hat{R}^j[\hat{I}_{q_1}^j, \dots, \hat{I}_{q_k}^j]$, we have at most $$m \cdot \left(\frac{n}{\hat{r}}+1\right)^k$$ such copies of $K_k$.
In addition, we could have more copies of $K_k$ not arising from this way. Let $K'$ be such a copy not yet counted. Then at least one edge $uv$ of $K'$ must be in $G^j[\widehat{V}_{j'}, \widehat{V}_{j''}]$ with $j'j''\notin \hat{R}^j$.
The reason why $j'j''\notin \hat{R}^j$ is one of two: either $G^j[\widehat{V}_{j'}, \widehat{V}_{j''}]$ has density less than $\hat{d}$, or it is not a regular pair.
The edges from the former cases extend to $\left(\left(1+\e^{1/2}\right)\frac{n}{r_j}\right)^{k-2}$ cliques, so this 
yields at most 
$$\sum_{i,i'\in [k] } \sum_{ j'\in \hat{I}^j_i, j''\in \hat{I}^j_{i'}} \hat{d}|\widehat{V}_{j'}||\widehat{V}_{j''}| \left(\left(1+\e^{1/2}\right)\frac{n}{r_j}\right)^{k-2}\leq 2\binom{k}{2} \hat{d} \left(\frac{n}{r_j}\right)^k$$
such copies of $K_k$ in $\hat{R}^j[\hat{I}_{1}^j, \hat{I}_{2}^j, \ldots, \hat{I}_{k}^j]$.
For the latter case, \ref{eq: common reg 5} yields that for fixed $j'$, there are at most $\eps_1 \hat{r}$  choices of $j''$ such that $G^j[\widehat{V}_{j'}, \widehat{V}_{j''}]$ is irregular.
So, by choosing $u$ which determines $j'$ and choosing $j''$ and $v\in \widehat{V}_{j''}$ and extending it to a copy of $K_k$, we can count that there are at most 
$$ \left(k\cdot \left(1+\e^{1/2}\right)\frac{n}{r_j}\right)\cdot \eps_1 \hat{r} \cdot \left(\frac{n}{\hat{r}}+1\right) \cdot \left(\left(1+\e^{1/2}\right)\frac{n}{r_j}\right)^{k-2}
\leq 2 k r_j \eps_1  \left(\frac{n}{r_j}\right)^{k}  $$
In sum, we have 
$$m \left(\frac{n}{\hat{r}}+1\right)^k + 2 \binom{k}{2} \hat{d} \left(\frac{n}{r_j}\right)^{k} + 2 k r_j \eps_1   \left(\frac{n}{r_j}\right)^{k}  \geq \left(\frac{d}{2}\right)^{\binom{k}{2}} \left(\frac{n}{r_j}\right)^k.$$
Since $\eps_1 \ll 1/r_j, \hat{d}$ and $\hat{d} \ll d, 1/k$, this implies that 
$$ m \geq \left(\frac{d}{3}\right)^{\binom{k}{2}} \left(\frac{\hat{r}}{r_j}\right)^{k} \geq 2k\hat{d} s_j^{k} .$$
This proves the claim.
\end{poc}

Using this claim, we can find a clique-walk transiting from $R^1$ to $\hat{R}^1$ and $\hat{R}^1$ to $R^2$. 
By our assumption, $[k]$ forms a clique in $R^1$.
For each $S\subseteq [k], \ell\in S, I_{\ell}\subseteq \hat{I}^1_{\ell}$, let $K(I_\ell: \ell\in S)$ be the number of copies of $K_k$ in $\hat{R}^1[\hat{I}^1_1,\dots, \hat{I}^1_k]$ which intersect $I_{\ell}$ for every $\ell\in S$.
When $I_{\ell} = \{i_{\ell}\}$, then we instead write $i_\ell$ to denote $\{i_{\ell}\}$.

By Claim~\ref{cl: clique counting} and the pigeonhole principle, there exists $i_k \in \hat{I}^1_k$ which belongs to at least $(k+1)\hat{d} s_1^{k-1}$ copies of $K_k$ in $\hat{R}^1[\hat{I}^1_1,\dots, \hat{I}^1_k]$. Now we choose vertices $i_k,i_1,\dots, i_{k-1}$ in order such that $i_\ell \in \hat{I}^1_{\ell}$ and the following holds for every $\ell\in [k-1]$.
\begin{enumerate}
    \item[(a)] $i_k, i_1,\dots, i_\ell$ lie in at least $(k-\ell+1) \hat{d} s_1^{k-1-\ell}$ copies of $K_k$ in $\hat{R}^1[\hat{I}^1_1,\dots, \hat{I}^1_k]$.
\end{enumerate}
Indeed, we can choose such $i_k,i_1,\dots, i_{k-1}$. By considering $i_0=i_k$, once we have $i_0,\dots, i_\ell$ for some $0\leq\ell \leq k-2$, we can choose next $i_{\ell+1}$ by simple pigeonhole property. Repeating this yields the desired $i_0,\dots, i_{k-1}$.
For each $\ell\in [k-1]$, let 
$$\hat{J}^{1}_{\ell} = \{ i\in \hat{I}^1_{\ell} : K(\{i_0,\dots, i_{\ell-1}, i\})\geq \hat{d} s_1^{k-1-\ell}\}.$$
Then we have 
\begin{align*}
    (k-\ell+2) \hat{d} s_1^{k-\ell} \leq  K(\{i_0,\dots, i_{\ell-1}\}) = \sum_{i\in \hat{I}^1_{\ell} } K(\{i_0,\dots, i_{\ell-1}, i\}) \leq |\hat{J}^{1}_{\ell}| s_1^{k-1-\ell} +  (s_1-|\hat{J}^{1}_{\ell}|)  \hat{d} s_1^{k-1-\ell}.
\end{align*}
 This yields that $|\hat{J}^{1}_{\ell}|\geq \hat{d} s_1$.

Now, we consider the sequence of the sets
$$U_1,U_2,\dots, U_{3k-1}  = V^1_1,\dots, V^1_{k}, \widehat{V}_{\hat{I}^1_1},\dots, \widehat{V}_{\hat{I}^1_{k-1}}, \widehat{V}_{i_k}, \widehat{V}_{i_1},\dots, \widehat{V}_{i_{k-1}}.$$
By Lemma~\ref{lem:tiny_part_good_enough}, for every distinct $i, j\in [2k-1]$ with $|i-j|\neq k$, the pair $(U_i,U_j)$ is an $(\eps_0^{1/2}, d+)$-regular pair in $G^1$ as $\frac{\eps_0}{2\hat{d}} \leq \eps_0^{1/2}$. 
For each $i\in [k+1,2k-1]$ and $j\in [2k,3k-1]$ with $j-i \neq  k$, the pair $(U_i,U_j)$ is an $(\eps_1^{1/3},\hat{d}+)$-regular pair in a graph $G'^{1}$ by Lemma~\ref{lem:union_regular} and \ref{eq: common reg 4}, for some subgraph $G'^{1}$ of $G^{1}$.
This sequence will be our starting pieces for the clique-walk.

We choose the ordered clique $(k,k-1,\dots, 1)$ of $R^2$. 
In the same way as before with the order of $(1,\dots,k)$ reversed, we can choose a clique $(i'_1,i'_k,\dots, i'_2)$ in $\hat{R}^2$ with $i'_j\in \hat{I}^2_j$ for each $j\in [k]$ and a sequence of the sets
$$U'_1,U'_2,\dots, U'_{3k-1}  = V^2_k,\dots, V^2_{1}, \widehat{V}_{\hat{I}^2_{j_k}},\dots, \widehat{V}_{\hat{I}^2_{j_{2}}}, \widehat{V}_{i'_1}, \widehat{V}_{i'_k},\dots, \widehat{V}_{i'_{2}}$$
in such a way that $(U'_i, U'_j)$ is $(\eps_0^{1/2},d+)$-regular pair in $G^2$ for $i\neq j\in [2k-1]$ with $|i-j|\neq k$. Additionally, for each $i\in [k+1,2k-1]$ and $j\in [2k,3k-1]$ with $j-i\neq k$, the pair $(U'_i,U'_j)$ is an $(\eps_1^{1/3}, \hat{d}+)$-regular pair in a graph $G'^{2}$ by Lemma~\ref{lem:union_regular} and \ref{eq: common reg 4}, for some subgraph $G'^{2}$ of $G^{2}$.

Applying Lemma~\ref{lem: clique_walk}, we obtain a clique-walk of length at most $3k^3$ in $\hat{R}^1$ from $(i_1,\dots, i_k)$ to $(i'_k,\dots, i'_1)$.
Let $c_1,\dots, c_{t}$ be such a sequence with $t\leq 3k^3$ is divisible by $k$. Then we define the sequence $U_1,\dots, U_{\ell}$ as follows:
$$\mathcal{U}=(U_1,\dots, U_{3k-1}, V_{i_k}, V_{c_1}, V_{c_2},\dots, V_{c_t}, V_{i'_1}, U'_{3k-1}, U'_{3k-2},\dots, U'_1, V^2_1,\dots, V^2_{k}, V^2_1,\dots, V^2_{k}).$$
In the last pieces, we repeat the segment $V^2_1,\dots, V^2_{k}$ until we obtain exactly $\ell$ sets in the family $\mathcal{U}$.

In this sequence $\mathcal{U}=(U_1,\dots, U_\ell)$ and the reduced graph $P^k_{[\ell]}$ on the vertex set $[\ell]$, for every $ij\in E(P^k_{[\ell]})$, the pair $(U_i,U_j)$ is $(\eps_0^{1/2},\hat{d}+)$-regular pair in either $G_1$ or $G_2$ or $G'_1$ or $G'_2$.

Let $J_1, J_2, J'_1, J'_2$ be the set of pairs $(i,j)$ with $i<j$ and $ij\in E(P^k_{[\ell]})$ where $(U_i, U_j)$ is $(\eps_0^{1/2},\hat{d}+)$-regular in $G_1, G_2, G'_1, G'_2$, respectively. In particular, we know that if $ij\in E(P^k_{[\ell]})$ satisfies 
$i\in [k+1,2k-1]$ and $ j\in [2k,3k-1]$, then $(i,j)\in J'_1$. If $i\in [t+3k+2,t+4k+1]$ and $j\in [t+4k+2,t+5k]$ then $(i,j)\in J'_2$. By remove some indices from $J_1, J_2$ if necessary, we may assume that $J_1,J_2, J'_1, J'_2$ forms a partition of $E(P^{k}_{[\ell]})$.

We encounter some technical problems here for applying the Lemma~\ref{lem:technical_blow}. In Lemma~\ref{lem:technical_blow}, two sets in $V_{g(1)},\dots, V_{g(\ell)}$ are either identical or disjoint, while some sets in $\mathcal{U}$ may not be.
Also, Lemma~\ref{lem:technical_blow} is stated for a single graph while the graphs $G'^j$ and $G^j$ may use common vertices, which renders some tedious technical problems.

In order to overcome this, we artificially make the sets in the sequence disjoint.
Take five disjoint subsets $Z_1, Z_2, Z_3, Z_4, Z_5$ of $V(G)$ with size $s_1\alpha' n, s_1\alpha'n, \alpha' n, s_2\alpha'n, s_2\alpha' n $ uniformly at random.  Lemma~\ref{lem:randomset_reg} yields that there is a choice of $Z_1,Z_2,Z_3,Z_4, Z_5$ such that for all $i,i'\in [5]$, the pair $(Z_i\cap U_j, Z_{i'}\cap U_{j'})$ is $(\eps_0^{1/30}, \hat{d}+)$-regular in $G^{a}$ if $(j,j')\in J_{a}$, and in $G'^{a}$ if $(j,j')\in J'_{a}$. 
Moreover, using Lemma~\ref{lem:concentration}, we can ensure that $|V^j_i\cap Z^a|$ and $|\hat{V}^j_i\cap Z^a|$ has size between $\frac{\hat{d} \alpha' n}{2\hat{r}}$ and $\frac{2\alpha' n}{\hat{d}\hat{r}}$.
We fix one such choice of $Z_1,Z_2,Z_3,Z_4, Z_5$.

Let $\mathcal{U}^*$ be the sequence $(U^*_1,\dots, U^*_\ell)$ as follows.
\begin{align}
U^*_i =\left\{
\begin{array}{ll}
U_i \cap Z_1 & \text{ if } i \in [k] \\
U_i \cap Z_2 & \text { if }i\in [k+1,2k-1] \\
U_i \cap Z_3 & \text { if }i\in [t+3k+2, t+4k+1] \\
U_i \cap Z_4 & \text { if }i\in [t+4k+2, t+5k] \\
U_i \cap Z_5 & \text { if }i\in [t+5k+1, \ell] \\
\end{array}
\right.
\end{align}
Let $G^*$ be the graph obtained from taking unions as follows.
$$G^*:=\bigcup_{a\in [2]} \bigcup_{(i,j)\in J_a} G^{a}[U^*_i, U^*_j]\cup \bigcup_{a\in [2]} \bigcup_{(i,j)\in J'_a} G'^{a}[U^*_i, U^*_j]. $$
Now, the sets in $\mathcal{U}^*$ are either disjoint or identical, and we have a single graph $G^*$ where $(U^*_i,U^*_j)$ is $(\eps_0^{1/30}, \hat{d}+)$-regular in $G^*$. Moreover, by the choices of $Z_1,\dots, Z_5$, we know that the size of each set $U^*_i$ is between $\frac{\hat{d} \alpha' n}{2\hat{r}}$ and $\frac{2\alpha' n}{\hat{d}\hat{r}}$.

From the assumption, $H$ admits the vertex partition $(P_{[\ell]}^k, \mathcal{W})$ and for each $i\in [\ell]$, we have $|W_i| \leq |H| \leq \ell \alpha n <  \frac{\alpha'n}{2\hat{r}} \leq |U^*_i|$ for each $i\in [\ell]$. 
We wish to apply Lemma~\ref{lem:technical_blow} to embed $H$ into $G^1\cup G^2$ on the regularity partition $( P^{k}_{[\ell]},\mathcal{U})$.
We have already checked two conditions \ref{blowup_assumption_0} and \ref{blowup_assumption_2}.
\ref{blowup_assumption_1} is trivial as $|V(H)|\leq 4\Delta \alpha n \leq  \frac{\alpha'n}{4\hat{r}} \leq \frac{1}{2}|U^*_i|$ for all $i\in [\ell]$.
\ref{blowup_assumption_3} and \ref{blowup_assumption_4} also follow from the assumptions of Connection Lemma as well as the definition of $\mathcal{U}$. Also the set $W_{[k]}$ which plays the role of $T$ in the application has size at most $k\alpha n$ with $k\alpha \ll \eps_1$. 
Hence, Lemma~\ref{lem:technical_blow} yields an embedding $\varphi: H \rightarrow G^1\cup G^2$ such that each $W_i$ is contained in $U_i$ and for each $x\in W_{[\ell]}$, we have $\varphi(x)\in A_x$ and for each $e\in E(\M)$, we have $|N_G(\phi(e))\cap U_{f(e)}| \geq \beta |U_{f(e)}|$.
\end{proof}

\section{Bandwidth theorem with target sets}
Assume that an $n$-vertex graph $H$ has an $\alpha n$-bandwidth partition $\mathcal{W}=(W_1,\dots, W_{t})$ of $H$. We say that a collection $\mathcal{I}=\{ I_i: i\in [p]\}$ is a \emph{collection of $\ell$-fragmented intervals} with respect to the partition $\mathcal{W}$ if each $I_i$ is a union of at most $\ell$ consecutive sets in $\mathcal{W}$ such that there are at least $\ell$ sets in $\mathcal{W}$ between any two intervals. 
 In other words,  we have $I_i = W_{[t_{2i-1},t_{2i}]}$ for each $i\in [p]$ and some sequence $t_1 \leq t_{2}\leq \dots \leq t_{2p}$ such that for each $i\in [p]$,
 $$|t_{2i}-t_{2i-1}|\leq \ell-1 \enspace \text{ and }\enspace |t_{2i+1}-t_{2i}|\geq \ell+1.$$
We say that this collection $\mathcal{I}$ is \emph{$\ell$-initial} if all intervals belong to the first $\ell$ sets in $\mathcal{W}$, i.e., we have $I_{[p]}\subseteq W_{[\ell]}$. 

Later in the applications (in particular, to execute Step 2-2 in the proof sketch), we will hand-pick some vertices $x$ of $H$ and embed them into $G$ in a specific way. This gives some restrictions on how we will later embed the neighbors of $x$. To reduce the restrictions, we will choose these vertices $x$ in such a way that the neighborhoods of $x$ belong to a collection of $\eta t$-fragmented intervals with respect to the bandwidth ordering for some small $\eta$. The advantage of each interval being short is that the restriction we get from this is not too strong. Also, the intervals being far apart ensures that the restrictions we get from one interval do not interfere with other intervals. This motivates the above definition of fragmented intervals. 

Developing the ideas in \cite{Boettcher2009}, Condon, Kim, K\"uhn and Osthus \cite{Condon2019} proved a bandwidth theorem for approximate decomposition. The proof of the following lemma is a modification of the proof in \cite{Condon2019}. As the proof in \cite{Condon2019} requires much more properties, our proof here is much simpler, but it still needs some additional care to ensure the desired properties. 

\begin{lemma}[Bandwidth theorem with target sets] \label{lem:technical_bandwidth}
Suppose $0 < 1/n \ll \a \ll \a' \ll \h \ll 1/r \ll \eps' \ll d'\ll d \ll \eps, 1/\D, 1/k, 1/s\leq 1$ and $k$ divides $r$. 
Let $H$ be an $n$-vertex graph with $\Delta(H)\leq \Delta$ such that $H$ has an $\alpha n$-bandwidth partition $\mathcal{W}=(W_1,\dots, W_t)$ with respect to an $\alpha n$-bandwidth ordering, and a proper $k$-coloring $c$ of $H$ (i.e., $\chi(H)\leq k$) such that $k$ divides $t$.

Let $G$ be an $n$-vertex graph and $\mathcal{V}=(V_1,\dots, V_r)$ is a partition of $V(G)$. Let $U$ be a vertex set not necessarily disjoint from $V(G)$. Assume we have a graph $G'$ on the vertex set $V(G)\cup U$ and $U_1,\dots, U_s$ is a partition of $U$ with $|U_i|\geq \eta n$ for each $i\in [s]$. 
Let $R$ be a graph on the vertex set $[r]$ with $\delta(R)\geq (1-\frac{1}{k}+\eps)r$ such that $K=(p^*_1,\dots, p^*_k)$ is an ordered clique in $R$. Let $\M$ be a multi-hypergraph on the vertex set $V(\M)=W_{[t-k+1,t]}$ together with a function $f:E(\M) \rightarrow [s]$.
Let $\mathcal{I}=\{ I_i: i\in [p]\}$ be a collection of $\eta t$-initial $\alpha' t$-fragmented intervals with respect to $\mathcal{W}$ such that $W_1\subseteq I_1$.
Suppose the following hold:
\stepcounter{propcounter}
\begin{enumerate}[label = {\bfseries \emph{\Alph{propcounter}\arabic{enumi}}}]
    \item $G$ admits $(\eps',d)$-regularity partition $(R,\mathcal{V})$ and $|V_i| = (1\pm \eps'^{1/2})\frac{n}{r}$ for each $i\in [r]$. \label{cond: G regularity}
    \item $\Delta(\M)\leq \Delta$ and the rank of $\M$ is at most $\Delta$. \label{cond: M hypergraph}
    \item For all $e\in \M$ and $x\in e$, if $x\in W_{t-k+i}$ for some $i\in [k]$ then $(V_{p^*_i},U_{f(e)})$ is an $(\eps',d+)$-regular pair in $G'$. \label{cond: M regularity}
    \item For each $I_j\in \mathcal{I}$ there exists an ordered $k$-clique $Q_j=(q_{j,1},\dots, q_{j,k}) \in \overrightarrow{K_k}(R)$ such that the following holds:  for each $x\in I_j$, there exists a set $A_x\subseteq V_{q_{j,c(x)}}$ of size at least $d'|V_{q_{j,c(x)}}|$. \label{cond: target sets}
\end{enumerate}
Then, there is an embedding $\p: H \rightarrow G$ such that the following hold:
\stepcounter{propcounter}
\begin{enumerate}[label = {\bfseries \emph{\Alph{propcounter}\arabic{enumi}}}]
    \item For all $I_j\in \mathcal{I}$ and $x \in I_j$, we have $\p(x) \in A_x$. \label{bandwidth_conclusion_1}
    \item For each $e \in E(\M)$, we have $|N_{G'}(\phi(e))\cap U_{f(e)}| \ge d'|U_{f(e)}|$. \label{bandwidth_conclusion_2}
\end{enumerate}
\end{lemma}

\begin{proof}[Proof of Lemma~\ref{lem:technical_bandwidth}]
By assumption, we have a hierarchy of constants as follows. 
\begin{align}\label{eq: hierarchy}
1/n \ll \a \ll \a' \ll \h \ll 1/r \ll \eps' \ll d' \ll d \ll \eps,1/\D,1/k,1/s \le 1.
\end{align}
The proof will go as follows. By slightly modifying the vertex partition $\mathcal{V}$, we make sure that $G$ admits a super-regularity partition $(Q,\mathcal{V})$ where $Q$ is a $K_k$-factor of $R$.
In Step 2, we will pick a number $t_1$ much smaller than $t$ and $\bigcup_{I_j\in \mathcal{I}}I_j \subseteq W_{[t_1]}$, and we will embed the first $t_1$ vertex sets $W_1,\dots, W_{t_1}$ in such a way that \ref{bandwidth_conclusion_1} holds. In Step 3, we will partition the remaining parts $W_{t_1+1},\dots, W_{t}$ of $H$ into $C_1, H_1, C_2, H_2,\dots, C_{\ell},H_{\ell}$ in such a way that each of $C_i, H_i$ contains consecutive $W_j$'s, and most of the sets $W_j$'s belong to $H_{[\ell]}= \bigcup_{i\in [\ell]} H_i$. The low bandwidth of $H$ together with the appropriate choice of sizes of $C_i$'s ensures that the graph $H_{[\ell]}$ has at least $\ell$ small components $H_1,\dots, H_{\ell}$. We plan to embed $H_{[\ell]}$ into $G$ utilizing the super-regularity structure $(Q,\mathcal{V})$. Since $Q$ is simply a disjoint union of $K_k$'s, $H_{[\ell]}$ having at least $\ell$ small components with $1/\ell \ll 1/r$ will be helpful for us. As we plan to use the blow-up lemma to embed $H_{[\ell]}$ into $G$ at the end of the proof, we ideally want to `distribute' the vertices in $H_{[\ell]}$ into $X_1,\dots, X_r$ so that $|X_i|=|V^*_i|$ holds, where $V^*_i$ is the left-over part of $V_i$ after embedding $W_{[t_1]}\cup C_{[\ell]}$. We first make a distribution of $H_{[\ell]}$ into $X^*_1,\dots, X^*_r$ so that $|X^*_i|= (1\pm \eta^{1/4})|V_i|$ in Step~3. 
In Step~4, we embed the parts $C_1,\dots, C_{\ell}$ utilizing the regularity structure $(R,\mathcal{V})$ of $G$. 
After Step~4, we have embedded all parts of $H$ but $H_{[\ell]}$, but the embedded parts are still negligible. One last issue here is that  $|X^*_i|$ and $|V^*_i|$ are not exactly equal. To ensure equality in their sizes, we want to move some vertices in $V_i$ to $V_j$ for some $i$ and $j$ while preserving the super-regularity structures in $Q$. However, these movements of vertices are only possible for specific pairs $(i,j)$. Thus, we will have to make some pre-arrangements at the beginning, which will be done in Step~1. In Step~5, we finally apply the blow-up lemma to embed the remaining graph $H_{[\ell]}$ into $G$ to obtain the final embedding.

Since $\delta(R -  V(K) )\geq (1-\frac{1}{k}+\eps)r -k \geq (1-\frac{1}{k})r$ and $k$ divides $r$, Theorem~\ref{thm:HS} implies that $R$ contains a $K_k$-factor containing $K$.
Let $Q^1,\dots, Q^{r'}$ be the vertex-disjoint copies of $K_k$ in $R$ and we give an arbitrary ordering of the vertices to assume that each $Q^i$ is an ordered clique and $Q^{r'}= K=(p^*_1,\dots, p^*_k)$.
We apply Lemma~\ref{lem: nice_regularity} to $G$, $R$ and $Q$ to get a new partition of $V(G)$.
By denoting this new partition as $\mathcal{V}$, we assume that the partition $\mathcal{V}=(V_1,\dots, V_r)$ of $V(G)$ satisfies the following.
\begin{equation}\label{eq: regularity Q}
    \begin{minipage}[c]{0.9\textwidth}
\begin{enumerate}
    \item[(a)] $|V_i| = \left(1\pm \eps'^{1/4} \right)n/r$ for all $i \in [r]$,
    \item[(b)] for each $ij \in E(Q)$, the graph $G[V_i,V_j]$ is $(\eps'^{1/4},d)$-super-regular, and
    \item[(c)] for each $ij \in E(Q)$, we have $||V_i|-|V_j||\le 1$.
    \item[(d)] $(R,\mathcal{V})$ is an $(\eps'^{1/4},d)$-regularity partition of $G$.
    \end{enumerate}
\end{minipage}
\end{equation}
This yields the desired super-regularity structure $(Q,\mathcal{V})$. In addition, the property (c) will be convenient for us later when we distribute the vertices of $H_{[\ell]}$ into $X^*_1,\dots, X^*_r$.

\mbox{}\\
\textbf{Step 1. Pre-arrangement for vertex redistribution on {\boldmath $\mathcal{V}$}.} 
Later, we will move some vertices across partitions $\mathcal{V}$ in $G$ while preserving the super-regularity structure $Q$. 
For each $j$, in order to be able to put a vertex $v$ into $V_j$, we need $|N_{G}(v)\cap V_{j'}|\ge (d-\e')|V_{j'}|$ for all $j'\in N_Q(j)$, as otherwise $(Q,\mathcal{V})$ is no longer a super-regularity partition for $G$ after putting $v$ into $V_{j}$. 
On the other hand, if $ij'\in E(R)$ for all $j'\in N_{Q}(j)$ for some $i\in [r]$, then Lemma~\ref{lem:large_super_regular} ensures that most of the vertices in $V_i$ can be added to $V_j$ while preserving the super-regularity structure of $Q$. 
So, we encode this information using the following auxiliary digraph $D$ on the vertex set $V(D)=[r]$ with the edge set \begin{equation}\label{eqn:directed}
E(D)=\left\{\vec{ij}:i\neq j\in [r], N_Q(j)\subseteq N_R(i)\right\}.
\end{equation}
In other words, $\vec{ij}\in E(D)$ implies that there is a large subset $\overline{V_i}\subseteq V_i$ such that any vertex in $\overline{V_i}$ can be moved to $V_j$ while somewhat preserving the super-regularity structure of $Q$. 
\begin{claim}\label{lem:directed}
There exists at least two $s\in[r']$ such that for any $j\in V(Q^s)$ and $i\in [r]$, there exists a directed path $P(i,j)$ from $i$ to $j$ in $D$.
\end{claim}
\begin{poc}
We first claim that any distinct $i,j\in [r]$ satisfy $N_D^+(i) \cap N_D^+(j)\neq \varnothing$. To see this, observe that $|N_R(\{i,j\})|\ge 2\d(R)-r\ge \left(1-\frac{2}{k}+\eps\right)r$. Thus, we have $$\big|\left\{a\in[r']:\left|N_R (\{i,j\})\cap V(Q^a)\right|\ge k-1\right\}\big|\ge \frac{\e r}{2}\ge 3.$$
Hence, there exists $a\in[r']$ such that $i,j\notin V(Q^a)$ while $\left|N_{R}(\{i,j\})\cap V(Q^a)\right|> k-2$. Finally, choose $j'\in V(Q^a)$ such that $Q^a\setminus \{j'\}\subseteq N_R(\{i,j\})$ and by \eqref{eqn:directed}, we have $j'\in N_D^+(i) \cap N_D^+(j)$, proving our claim.

Let $A(j)=\left\{i\in[r]:\text{there exists a directed path from   } i \text{  to  } j \text{  in  } D\right\}$. Fix a $j^*\in[r]$ that maximizes $|A(j^*)|$. Suppose there exists $i\notin A(j)$. Then, by the previous claim, we know that there exists $j'\in N_D^+(i) \cap N_D^+(j^*)$. Thus, $A(j^*)\cup\{i\}\subseteq A(j')$, contradicting the maximality of $|A(j^*)|$. Hence, we must have $A(j^*)=[r]$.

Since $\d(R)\ge \left(1-\frac{1}{k}+\frac{\eps}{2}\right)r$, we have $\left|\left\{a\in[r']:\left|N_{R}(j^*)\cap V(Q^a)\right|= k\right\}\right|\ge \frac{\eps r}{2} \ge 2$. Thus, there exist at least two choices of $s \in [r']$ such that $V(Q^s)\subseteq N_R(j)$. This implies that $V(Q^s)\subseteq N_D^+(j^*)$. Finally, for every $j\in V(Q^s)$ and $i\in[r]$, there exists a directed path from $i$ to $j^*$ in $D$ and a directed edge from $j^*$ to $j$, thus there is a directed path from $i$ to $j$ in $D$. This finishes the proof of Claim~\ref{lem:directed}.
\end{poc}
By re-indexing if necessary, we assume that $Q^1 =(1,2,\dots, k)$ satisfies Claim~\ref{lem:directed} to assume the following. Note that Claim~\ref{lem:directed} yields more than one choice of $s$, so we can assume that $Q^1$ is different from the clique $Q^{r'}=K$.

\begin{align}\label{eq: dipath in D}
\text{For any $j\in [k]$ and $i\in [k+1,r]$, there exists a directed path $P(i,j)$ from $i$ to $j$ in $D$.}
\end{align}

\mbox{}\\
\textbf{Step 2. Embedding the first segments of vertices.} 
The vertices in $I_j \in \mathcal{I}$ have to be embedded into $V_{Q_j} = \bigcup_{q\in Q_j} V_q$. However, $Q_j$ and $Q_{j+1}$ can be different cliques of $R$ that are far apart. Since $\mathcal{I}$ is $\alpha' t$-fragmented, there are many vertices between the interval $I_j$ and $I_{j+1}$. Lemma~\ref{lem: clique_walk} provides a clique-walk $P_j$ between $Q_j=(q_{j,1},\dots, q_{j,k})$ and $Q_{j+1}=(q_{j+1,1},\dots, q_{j+1,k})$ and we can utilize these clique-walks to embed the vertices between two intervals, establishing the connection. We will achieve this while embedding the first segments $W_1,\dots, W_{t_1}$ of $H$ for some $t_1$ with $t_1 \leq  (\eta+2k\alpha ) t$.
For our convenience, let $Q_{p+1} = (1,\dots,k)$ be the ordered clique.

We choose a minimum $t_1$ so that
$W_{[t_1]}=\bigcup_{i\in [t_1]} W_i$ contains every intervals in $\mathcal{I}$ and $k$ divides $t_1$. As $\mathcal{W}$ is an $\alpha n$-bandwidth partition of $H$, the set $W_{[t_1]}$ may contain at most $(\eta+2k\alpha) n$ vertices.
In this step, we embed the vertices of $W_{[t_1]}$ in $G$.

Let $\ell_1=0$ and $\ell'_1$ be the minimum number divisible by $k$ so that 
$I_1\subseteq W_{[\ell'_1-k]}$.
For each $i\in [p]\setminus\{1\}$, choose the maximum $\ell_i$ and the minimum $\ell'_i$ so that $I_i\subseteq W_{[\ell_i+k+1, \ell_{i}+\ell'_i-k]}$ and both $\ell_i,\ell'_i$ are divisible by $k$.
Since $\mathcal{W}$ is an $\alpha n$-bandwidth partition of $H$, we have $\alpha^{-1} \leq t \leq 2k \alpha^{-1}$.
Hence, the facts that $0<\alpha \ll \alpha' \ll 1/k \leq 1$ and that $\mathcal{I}$ is $\alpha' t$-fragmented imply the following.
\begin{enumerate}
	\item[(a)] $\ell_1< \dots < \ell_p$ and
     \item[(b)] $\ell'_i \leq 2\alpha' t$ and $\ell_{i+1}-(\ell_i+\ell'_i) \geq \alpha't/2$.
\end{enumerate}
We now collect the indices between $\ell_i$ and $\ell_i+\ell_i'$. For each $i\in [p]$, we define 
$$J_i = \{ j\in [t_1]: \ell_{i}+1 \leq j \leq \ell_{i}+\ell'_1\}.$$
Let $\ell_{p+1} = t_1+k$.
We next collect the indices between $J_i$ and $J_{i+1}$ to obtain the following set $J'_i$. For each $i\in [p]$, we define
$$J'_i = \{ j\in [t_1]: \ell_{i}+\ell'_i+1 \leq j \leq \ell_{i+1}\}.$$
Note that for each $i>1$, if $j$ is one of the smallest $k$ indices or largest $k$ indices in $J_i$, then $W_j$ does not intersect with $I_i\in \mathcal{I}$. These buffer indices will later make sure the `moreover' part of \ref{blowup_assumption_4} in our application of the blow-up lemma (first $k$ sets in $\{W_{j}: j\in J_i\}$ and the vertices in $I_i$ will play the role of $T$ and the last $k$ sets in $\{W_{j}: j\in J'_i\}$ whose vertices may have neighbors in the first $k$ sets in $\{W_{j}: j\in J_{i+1}\}$ will play the role of the vertices belonging to the edges of $\M$).

 We wish to embed the vertices in $W_{J_i} = \bigcup_{j\in J_i} W_j$ into $V_{Q_i}=\bigcup_{q\in Q_i} V_q$ while embedding the vertices in $W_{J'_i}=\bigcup_{j\in J'_i} W_j$ into $W_P=\bigcup_{q\in P} W_q$ where $P$ is a clique-walk from $Q_i$ to $Q_{i+1}$ which we will later obtain using Lemma~\ref{lem: clique_walk}.

If $m>1$, for each $j\in [k]$ and each vertex $x\in W_{\ell_{m}+j}$ in the first $j$-th set to be embedded next, we let $N_{m,x} = N_H(x)\cap \bigcup_{i\leq \ell_{m}} W_i$ and let $A_x:= N_G(\phi_{m-1}( N_{m,x}))\cap V_{q_{m,j}}$.

Assume that we have found an embedding $\phi_{m-1}$ of all vertices in $W_1,\dots, W_{\ell_{m}}$ into $G$, and additionally, assume the following.
\stepcounter{propcounter}
\begin{enumerate}[label = {\bfseries \emph{\Alph{propcounter}\arabic{enumi}}}]
	\item The vertices $x\in \bigcup_{i\in [m-1]} I_i$ are embedded in $A_x$ (as in the assumption of Lemma~\ref{lem:technical_bandwidth}).\label{cond: assmp1} 
	\item If $m>1$, then the sets $W_{\ell_{m}-k+1}, W_{\ell_{m}-k+2}, \dots, W_{\ell_{m}}$ are embedded in $V_{q_{m,1}},\dots, V_{q_{m,k}}$, respectively, where $Q_m=(q_{m,1},\dots, q_{m,k})$.\label{cond: assmp2}
	\item  For each $j\in [k]$ and each vertex $x\in W_{\ell_{m}+j}$, we have $|A_x| \geq \frac{1}{2}d'|V_{q_{m,j}}|.$ \label{cond: assmp3} 
\end{enumerate}
Note that the set $W_{[\ell_{m}+1,\ell_{m}+k]}$ is disjoint from the intervals in $\mathcal{I}$ by the definition of $\ell_i$ and $\ell'_i$, hence the set $A_x$ is uniquely defined for the vertices $x\in W_{[l_m+1,l_m+k]}\cup I_m$. 
Note that an empty embedding $\phi_{0}$ vacuously satisfies all three conditions above as
\ref{cond: target sets} implies  \ref{cond: assmp3} when $m=1$.

Now we will find an embedding $\phi_{m}$ extending $\phi_{m-1}$ which also satisfies the above three conditions. 
For each $j\in J_m$, we assign $W_j$ to $X'_{q_{m,j'}}$ where $j = j'~({\rm mod}~k)$ with $j'\in [k]$. As $k$ divides both $\ell_m$ and $\ell'_m$, the last $k$ sets $W_{\ell_m+\ell'_m-k+1},\dots, W_{\ell_m+\ell'_m}$ are assigned to $X'_{q_{m,1}},\dots, X'_{q_{m,k}}$, respectively.

By applying  Lemma~\ref{lem: clique_walk} with the ordered clique $Q_m$ and $Q_{m+1}=(q_{m+1,1},\dots, q_{m+1,k})$, we obtain a clique-walk $P$ satisfying the following.
\begin{enumerate}
	\item $P = p_{1}\dots p_{a}$ with $3k\leq a \leq 3k^3$ and $k$ divides $a$.
	\item for all $i,j\in [a]$ with $|i-j|\leq k-1$, we have $p_{i}p_{j}\in E(R)$.
	\item for each $i\in [k]$, $p_{i}= r_{i}$ and $p_{a-i+1} = q'_{k-i+1}$.
\end{enumerate}
Now for each $j\in J'_m$, we assign the set $W_{j}$ to $X'_{j'}$ for some $j'$.
We assign each of the sets
$$W_{\ell_{m}+\ell'_m+1}, W_{\ell_{m}+ \ell'_m+2},\dots, W_{\ell_{m+1}}$$ 
in order one by one to the vertex sets $X'_{p_1}, X'_{p_2},\dots, X'_{p_{a-1}}, X'_{p_{a}}, X'_{p_{a-k+1}}, X'_{p_{a-k+2}},\dots$ respectively. Note that once we reach $X'_{a}$, we repeat the segment $X'_{a-k+1},\dots, X'_{a}$. As $\ell_{m+1}-(\ell_m+\ell'_m)$ is a number bigger than $\frac{\alpha'}{2\alpha} > a$ which is divisible by $k$ and $a$ is also divisible by $k$, we know that the last $k$ sets $W_{\ell_{m+1}-k+1},\dots, W_{\ell_{m+1}}$ are assigned to $X'_{q_{m+1,1}},\dots, X'_{q_{m+1,k}}$, respectively in order.
Note that $p_i=p_j$ is possible for distinct $i,j$, so each $X'_j$ is possibly a union of several $W_j$s. Let $J$ be the set of indices $j$ such that $X'_j$ is not empty, then $|J|\leq a\leq 3k^3$.

Let $W = W_{J_m\cup J'_m}$. 
Now consider induced subgraph $H'=H[W]$ of $H$ with the vertex partition $(X'_j )_{j\in J}$. Consider the vertex partition $(V_j)_{j\in J}$ of $G[V_{J}]$ where $V_{J}= \bigcup_{j\in J} V_j \setminus \phi_{m-1}(W_{[\ell_m]})$. 

For each $j\in [k]$ and a vertex $x \in W_{\ell_{m+1}+j}$, we add a hyperedge $e= N_{m,x}$ to the multi-hypergraph $\M'$ and let $f(e)= j$ and let $U_{j}= V_{q_j}$.
Let $R'$ be the graph with the vertex set $V(P)$ and the edge set $E(P)$.
Now we apply Lemma~\ref{lem:technical_blow} with the following roles of parameters.
{\begin{center}\small
\begin{tabular}{ |c|c|c|c|c|c|c|c|c|c|c|c|c|c| }
\hline
 objects/parameters & $G[V_J]$ & $H'$ & $R'$ & $\M'$ & $U_j$ & $I_{m}\cup W_{[\ell_{m}+1,\ell_{m}+k]}$ & $A_x$ & $\eps'^{1/2}$ & $d'/3$ & $d'$ & $a$ & $\eta^{1/2}$ \\ 
 \hline
 playing the role of & $G$ & $H'$ & $R$ & $\M$ & $U_j$ & $T$ & $A_x$ & $\eps$ & $\alpha$ & $\beta$ & $k$ & $\gamma$ \\ \hline
\end{tabular}
\end{center}}
Note that $|W_{[\ell_m]}|\leq (\eta+2k\alpha) n$ and $\eta \ll 1/r$. Hence the set $|\phi_{m-1}(W_{[\ell_m]})|$ is much smaller than one set $V_j$, so we have $|X'_j| < \eta^{1/2} |V_j|$. 
We next check that the application of Lemma~\ref{lem:technical_blow} is possible. Indeed, when $g$ is naturally chosen to make correspondence with $X'_{p_j}$ to $V_{p_j}$, \ref{blowup_assumption_1} holds. 
Furthermore, \ref{blowup_assumption_0} holds by Lemma~\ref{lem:tiny_part_good_enough} and by the properties of the clique-walk $P$.  The maximum degree condition on $H$ together with the definition of $\M'$ ensures \ref{blowup_assumption_2}.
 The definitions of $U_j$ and $f$ ensure \ref{blowup_assumption_3} and the assumption on $I_{m}$ and \ref{cond: assmp3} yield \ref{blowup_assumption_4}.
 Then we obtain an embedding $\phi'$ of $H'$ into $G[V_J]$ such that the followings hold.
 \begin{enumerate}
 	\item[(a)] for all $x\in I_{m} \cup \bigcup_{j\in [k]} W_{\ell_{m}+j}$, we have $x\in A_x$.
 	\item[(b)] for each $e\in E(\M')$, we have $|N_G(\phi'(e))\cap U_{f(e)}| \geq d'|U_{f(e)}|$.
 \end{enumerate}
 The first property for $x\in \bigcup_{j\in [k]} W_{\ell_{m}+j}$ ensures that the union $\phi_{m} = \phi_{m-1}\cup \phi'$ is indeed an embedding of the first parts $H[W_{[\ell_{m+1}]}]$ into $G$.
 The first property for $x\in I_{m}$ implies that \ref{cond: assmp1} holds for $\phi_{m}$. By our definition of $X'_{j}$, \ref{cond: assmp2} holds for $\phi_{m}$. Finally, the second property above together with the definition of $f$ ensures that \ref{cond: assmp3} holds for $\phi_{m}$. 
 Hence, we may repeat this until we obtain $\phi_{p}$ embedding all vertices in $W_1,\dots, W_{t_1}$ while \ref{cond: assmp1}--\ref{cond: assmp3} hold.
 In particular, our assumption $Q_{p+1}=(1,2,\dots, k)$ implies the following.
\begin{equation}
\begin{minipage}[c]{0.9\textwidth}
	For each $j\in [k]$ and each vertex $x\in W_{t_1+j}$, consider the set $N_{m,x} = N_H(x)\cap \bigcup_{i\leq t_1} W_i$, then we have
	$|N_G(\phi_{p}( N_{p,x}))\cap V_{j}|\geq d'|V_{j}|.$
	\label{eq: t_1 end}
\end{minipage}	
\end{equation}
In addition, since $k$ divides $t_1$, the definition of $\alpha n$-bandwidth partition ensures that the vertices in $W_{t_1+j}$ are colored $j$ by the proper coloring $c$ for each $j\in [k]$.

\mbox{}\\
\textbf{Step 3. Distributions of the remaining vertices of {\boldmath $H$}.} 
Let $H_0 = H[W_{[t_1]}]$ be the subgraph of $H$ embedded in Step~2.
For each $i\in [r]$, let $V'_i = V_i\setminus \phi_p(V(H_0))$ be the set of unused vertices in $V_i$.
Since at most $|V(H_0)|\leq (\eta+2k\alpha)n< \frac{\eps'^{1/2} n}{r}$ vertices have been embedded, by \eqref{eq: regularity Q}, the number $n'_i:=|V'_i|$ of vertices in each set $V'_i$ satisfies that for each $i\in [r]$,
\begin{align}\label{eq: size n'}
n'_i=\left(1\pm 2\e'^{1/4}\right)\frac{n}{r} \enspace \text{ and } \enspace |n'_i - n'_j|\leq 2\eta n \text{ for $ij\in E(Q)$}.
\end{align}
Let $W=\bigcup_{i>t_1} W_i$ be the set of unembedded vertices and let $H'=H[W]$.
We again partition the set $[t_1 +1,t]$ into 
intervals $L_{1,1},L_{1,2}, L_{1,3}, L_{1,4}, L_{2,1},\dots, L_{\ell,1},L_{\ell,2},L_{\ell,3}, L_{\ell,4}$ in order satisfying the following.
\begin{enumerate}
\item[(a)] $\alpha' n \leq |W_{L_{i,j}}| \leq 2\alpha' n$ for $i\in [\ell]$ and $j\in [3]$
\item[(b)] $2\eta n \leq |W_{L_{i,3}\cup L_{i,4}\cup L_{i+1,1}}| \leq 3\eta n$ for each $i\in [\ell]$, where $L_{\ell+1,1}=\varnothing$.
\item[(c)] $|L_{i,j}|$ is a multiple of $k$ for all $i,j$. 
\end{enumerate}
Since each $W_i$ has size at most $k\alpha n$ and $0<\alpha \ll \alpha'\ll \eta<1$, it is easy to see that such intervals exist. For each $i\in [\ell]$, we define the following sets and graphs: 
$$T_0= W_{L_{1,1}}, \enspace T_i= W_{L_{i+1,1}}, \enspace C_i=H[W_{L_{i,2}}], \enspace T'_i=W_{L_{i,3}} \enspace \text{ and }\enspace H_i = H[W_{L_{i,3}}\cup W_{L_{i,4}}\cup W_{L_{i+1,1}}].$$
Note that these choices yield $\frac{1}{4\eta} \leq \ell \leq \frac{1}{\eta}.$
Furthermore, the set $V(\M)$ belongs to $L_{\ell,4}$, hence belongs to $H_{\ell}$ and is disjoint from $T'_{\ell}$.

Then most of the vertices in $H$ belong to the graph $H_{[\ell]}=\bigcup_{i\in [\ell]} H_i$, as the other parts $C_i$ have at most $2\alpha'n \ell \leq \alpha'^{1/2} n$ vertices in total. In addition, the graph $\bigcup_{i\in [\ell]} H_i$ has at least $\ell$ connected components of size at most $3\eta n$ because $|C_i|\ge \alpha'n$ while the bandwidth of $H$ is smaller than $\alpha n< \alpha' n$, implying that there are no edges between $H_i$ and $H_j$ with $i\neq j$. 

Recall that $H$ has exactly the same number of vertices as $G$. In order to use the blow-up lemma (Lemma~\ref{lem:technical_blow}) to embed the graph $H_{[\ell]}$, we want to use the fact that $G$ admits the $\eps$-super-regular partition $(Q,\mathcal{V})$. 
In order to do this, we again want to embed $C_i$ to establish the connections between disjoint cliques in $Q$ utilizing the clique-walks in $R$ as before. 
The parts $C_i$ are the parts that we will embed using the clique-walks in $R$.
The sets $T_i$ and $T'_i$ are the vertices $x$ which potentially has neighbors in $C_i$, so we have to embed it to a target set $A_x$ which is a common neighborhood of the embedded images of the neighbors of $x$ in $C_i$. 
In order to embed $H_{[\ell]}$ into $G$, we need to first assign the vertices in $H_i$ into the cliques in $Q$ so that roughly $n'_j$ vertices are assigned to be embedded into $V'_j$. We will make such a distribution now.

Partition the vertices in $H_i$ (similarly $T_i,C_i,T'_i$) into $H_{i,1},\ldots,H_{i,k}$ where the vertices in $H_{i,j}$ are colored $j$ with the coloring $c$. By the definition of $\mathcal{W}$, we know that for any $j$ each of $k$ consecutive sets $W_j, W_{j+1},\dots, W_{j+k-1}$ within $H_{i}$ belongs to different parts of $H_{i,1},\ldots,H_{i,k}$.

Note that our previous embedding $\phi_{p}$ satisfies \eqref{eq: t_1 end}.
For each $x\in T_0$ with $x\in W_{t_1+j}$ and $j\in [k]$, let 
$A_x = N_G(\phi_{p}( N_{p,x}))\cap V_{j}$. In particular, if $x$ has no neighbors in the already embedded subgraph $H_0$ of $H$, then $A_x = V_{j}$. As $d'<1$, for all $x\in T_0$ with $c(x)=j$, we have
$$|A_x|\geq d'|V_j|.$$ 
Recall that $r'=\frac{r}{k}$ and $Q^1,\dots, Q^{r'}$ are the disjoint ordered $k$-cliques in $Q$ and $Q^1= (1,2,\dots, k)$. We will distribute the sets $W_j$ within $H_{[\ell]}$ to obtain the sets $X^*_1,\dots, X^*_r$ so that the set $X^*_i$ is slightly bigger than $n'_i=|V'_i|$ if $i\in [k]$, and slightly smaller than $n'_i$ if $i\in [k+1,r]$.
 The property \eqref{eq: dipath in D} will later ensure that we can move some vertices in $V'_j$ to $V'_i$ with $j\in [k+1,r]$ and $i\in [k]$ so that we can later obtain the exact equality in \ref{blowup_assumption_1} for the final application of the blow-up lemma~(Lemma~\ref{lem:technical_blow}).

\begin{claim}\label{cl:rearrangement_of_H}
There exists a partition $\mathcal{X}^*= (X^*_1,\dots, X^*_r)$ of the set $\bigcup_{i\in[\ell]} V(H_i)$ such that the followings hold.
\begin{enumerate}[leftmargin=*]
    \item[(a)] The graph $\bigcup_{i\in [\ell]}H_i$ admits the vertex partition $(Q,\mathcal{X}^*)$.
    \item[(b)] Each $X^*_i$ is a union $\bigcup_{j\in S} H_{j,\ell_j}$ of sets for some $S\subseteq [\ell]$ and some $\ell_j \in [k]$ for each $j\in S$.
    \item[(c)] For each $i\in [r]\setminus [k]$, we have 
    $$(1-\eta^{1/4})n'_i \le |X_i| \le (1-\h^{2/3})n'_i.$$
    \item[(d)] For each $i\in [k]$, we have 
    $$(1+\eta^{2/3}) n'_i \leq |X_i| \leq (1+\eta^{1/4})n'_i.$$
    \item[(e)] For each $i\in [k]$, the set $W_{t-k+i} \subseteq X^*_{p^*_i}$.
\end{enumerate}
\end{claim}
\begin{poc}
We assign $H_{\ell, i}$ to $X_{q_{r',i}}=X_{p^*_i}$ for each $i\in [k]$, and let $\pi_{\ell+1}$ be the identity permutation on $[k]$.
For each $i\in[\ell-1]$, we independently at random choose a number $s'\in[r']$ such that $s'\in[r']$ is chosen with probability 
\begin{align*}
\begin{cases*}
p_{s'}= (1-\h^{1/2})\frac{1}{n}\sum_{i\in V(Q^{s'})} \n_i & for $s'>1$,\\
1-\sum_{s'>1} p_{s'} & for $s'= 1$
\end{cases*}
\end{align*}
and choose a permutation $\pi_s:[k]\rightarrow [k]$ uniformly at random. Then, for each $j\in[k]$, we add $H_{i,j}$ to $X_{q_{s,\pi_s(j)}}$ where $q_{s,\pi_s(j)}$ is the $\pi_s(j)$th vertex in the ordered clique $Q^s$.

The first two properties are immediate from the construction and the remaining two properties hold with high probability by a standard application of the Azuma's inequality together with \eqref{eq: size n'}. The last property follows from how we distribute $H_{\ell,i}$. 
This proves the claim. 
\end{poc}
Let $\mathcal{X}^* = (X^{*}_1,\dots, X^{*}_r)$ be the partition of $\bigcup_{i\in[\ell]} V(H_i)$ obtained from Claim~\ref{cl:rearrangement_of_H}.

\mbox{}\\
\textbf{Step 4. Embedding {\boldmath $C_i$} providing connections between {\boldmath $V_{Q^i}$}.} 
For each $i\in [\ell]$, assume that each $H_{i,1},\dots, H_{i,k}$ are assigned into $X_{q_{s_i, \pi_i(1)} },\dots, X_{q_{s_i,\pi_i(k)}}$ in Claim~\ref{cl:rearrangement_of_H} where $(q_{s_i,\pi_i(1)},\dots, q_{s_i, \pi_i(k)})$ is an ordered $k$-clique which we call $Q(i)$ with $V(Q(i)) = V(Q^{s_i})$.  Let $Q(0)=Q^1=(1,2,\dots, k)$, so $q_{s_0,\pi_0(j)} = j$ for each $j\in [k]$.
Now we will embed $T_0 \cup C_1$ using a clique-walk from $Q(0)$ to $Q(1)$ and embed $C_i$ using a clique-walk from $Q(i-1)$ to $Q(i)$ for each $i\geq 2$.

For each $i\in [\ell]$, apply Lemma~\ref{lem: clique_walk} on $R$ to find a clique-walk $P^i$ satisfying the following.
\begin{enumerate}
	\item[(a)] $P^i = p^i_{1}\dots p^i_{a_i}$ with $3k\leq a_i \leq 3k^3$ and $k$ divides $a_i$.
	\item[(b)] for all $j,j'\in [a_i]$ with $|j-j'|\leq k-1$, we have $p^i_{j}p^i_{j'}\in E(R)$.
	\item[(c)] for each $j\in [k]$, $p^i_{j}= q_{s_{i-1},\pi_{i-1}(j)}$ and $p^i_{a_i-k+j} = q_{s_{i}, \pi_i(j) }$.
\end{enumerate}
Let $L_{i+1,2}$ be $\{g+1,\dots, g+\ell'\}$ for some $g,\ell'$, then as $t_1$ is a multiple of $k$ and $|L_{i,j}|$ is a multiple of $k$ for all $i,j$, both $g$ and $\ell'$ are multiples of $k$.
We assign the sets $W_{g+1},\dots, W_{g+\ell'}$ to the sets $X'_{p^{i}_1}, X'_{p^{i}_2},\dots, X'_{p^{i}_{a_i}}, X'_{p^{i}_{a_i-k+1}}, X'_{p^{i}_{a_i-k+2}},\dots, X'_{p^{i}_{a_i}}$ one by one in this order. Note that once we reach the set $X'_{p^{i}_{a_i}}$, then we repeat $X'_{p^{i}_{a_i-k+1}}, X'_{p^{i}_{a_i-k+2}},\dots, X'_{p^{i}_{a_i}}$ again. As $\ell'$ and $a_i$ are both divisible by $k$, the last set $W_{g+\ell'}$ is assigned to $X'_{p^{i}_{a_i}}$.

We now aim to embed each $C_i$ into $G$ so that the vertices $x$ in $T_i\cup T'_i$ can be later embedded so that the edges between $C_i$ and $T_i\cup T'_i$ are also embedded into $G$. For this, we want the neighbor set $N_H(x)\cap C_i$ to be mapped into $G$ so that their images have a large common neighborhood in $G$.
Consider the hypergraph $\M^i$ with the edge set $\{N_{H}(x)\cap V(C_i) : x\in T_i\cup T'_i\}$. 
For each $j\in [k-1]$ and $j'\in [2k-2]\setminus [k-1]$, we let 
$$U^i_j = V_{p^i_j} \text{ and } U^i_{j'} =V_{p^i_{a_i-k+1+j'}}.$$ 
Now we define function $f^i$ on $E(\M^i)$. 
For each $j\in [k-1]$ and $j'\in [2k-2]\setminus [k-1]$ and the vertices $x\in W_{g-k+j}$ and $y\in W_{g+\ell'+j'-k+1}$, if $e=N_H(x)\cap V(C_i)$ then we let $f^i(e)= j$ and if $e= N_H(y)\cap V(C_i)$ then we let $f^i(e)=j'$. Note that multi-edges $e,e'$ with the same vertex sets could have different $f^i$-values based on where they come from.

Let $\phi'_0=\phi_{p}$. Assume that for some $0\leq m<\ell$, we have $\phi'_{m}$ which extends $\phi_{p}$ satisfying the following. 
\stepcounter{propcounter}
\begin{enumerate}[label = {\bfseries \emph{\Alph{propcounter}\arabic{enumi}}}]
\item If $m>0$, then $\phi'_m$ embeds $H_0\cup C_{[m-1]} \cup T_0$ into $G$.
    \item For each $i< m$ and an edge $e\in \M^i$, $|N_G(\phi'_{m}(e))\cap U^i_{f(e)}|\geq d' |U^i_{f(e)}|$. \label{enough_neighbor_C}
\end{enumerate}
Note that this embedding trivially exists when $m=0$.
Let $W=\varnothing$ if $m=0$, otherwise let $W= C_{[m-1]} \cup T_0$.
For each $i\in [r]$, let $V^*_i = V'_i\setminus \phi'_{m}(W)$ be the set of unused vertices in $V'_i$ and let $V^*=\bigcup_{i\in [r]}V_i^*$. Since $\bigcup_{i<m} C_i$ contains at most $\alpha'^{1/2}n <  \frac{\eta n}{r}$ vertices, we have $|V^*_i| = (1\pm 3\eps'^{1/2})\frac{n}{r}$.

In order to embed $C_{m+1}$, we apply Lemma~\ref{lem:technical_blow} with the following roles of parameters.
\begin{center}
\begin{tabular}{ |c|c|c|c|c|c|c|c|c| }
\hline
 sets/graphs & $G[V_{V(P^{m+1})}]$ & $H[C_{m+1}]$ & $V^*_{j}$ & $U^{m+1}_j$  & $X'_j$ & $\substack{T_0 \text{ if $m=0$ }\\ \varnothing~ \text{ otherwise} }$ & $A_x \cap V^*$ \\ 
 \hline
 playing the role of & $G$ & $H$ & $V_j$ & $U_j$  & $X_j$ & $T$ & $A_x$ \\ \hline
\end{tabular}
\end{center}
Below we check that the hypotheses of Lemma~\ref{lem:technical_blow} hold with the following parameters in all of the above cases.
\begin{center}
\begin{tabular}{ |c|c|c|c|c|c|c|c|c|c| }
\hline
 parameters/numbers & $2n/r$ & $\a'^{1/2}r$ & $\e'$ & $d'/2$ & $d'$ & $d$ & $\D$ & $k$ & $2k$ \\ 
 \hline
 playing the role of & $n$ & $\g$ & $\e$ & $\a$ & $\b$ & $d$ & $\D$ & $k$ & $s$ \\ \hline
\end{tabular}
\end{center}
Conditions \ref{blowup_assumption_0} and \ref{blowup_assumption_1} hold by a simple application of Lemma~\ref{lem:tiny_part_good_enough} and the natural choice of the function $g$. Conditions \ref{blowup_assumption_2}-\ref{blowup_assumption_3} are easy to check. If $m=0$, then Condition \ref{blowup_assumption_4} follows from \eqref{eq: t_1 end}, otherwise it follows from \ref{enough_neighbor_C} and the the fact that $|A_x\cap V^*|\geq |A_x| - \alpha'^{1/2}n \geq \frac{1}{2}|A_x|$. Thus, we can find an embed $\phi'_{m+1}$ of $H[C_i]$ into $G[V^*]$ satisfying \ref{enough_neighbor_C}. By repeating this, we can embed all $T_0,C_1,\dots, C_{\ell}$. Define $T=\cup_{i\in [l]} (T_i \cup T'_{i-1})$. Finally, the last discussion shows that after the embedding of all the vertices in $T_0\cup C_{[\ell]}$, we have the following.
\begin{equation}
	\begin{minipage}[c]{0.9\textwidth}
	For every $i\in [\ell]$ and $x\in (T'_{i-1}\cup T_{i})$ which is colored $j$ by the coloring $c$, the set ${A_x=N_G(\phi'_{\ell}(N_H(x)\cap C_i))\cap V_{q_{s_i,\pi_i(j)}}}$ satisfies $|A_x|\geq d' |V_{q_{s_i,\pi_i(j)}}|$. \label{available_vertex_C}
\end{minipage}
\end{equation}

\mbox{}\\
\textbf{Step 5. Rearrangements of some vertices in {\boldmath $\mathcal{V}$}.} 
Let $H^*$ be the graph $H[ V(H_0)\cup T_0\cup C_{[\ell]}]$ embedded so far.
For each $i\in [r]$, let 
$$V^*_i =  V'_i \setminus \phi'_{\ell}( T_0\cup C_{[\ell]}) =V_i \setminus \phi'_{\ell}( V(H^*)).$$
Letting $n^*_i = |V^*_i|$, as $|T_0\cup C_{[\ell]}|\leq \alpha'^{1/2}n$, we have
\begin{align}\label{eq: size n^*}
n^*_i= n'_i \pm \alpha'^{1/2} n = (1\pm \alpha'^{1/3}) n'_i.
\end{align}
Now the vertices in $H_0$ and $C_i$ are all embedded by $\phi'_{\ell}$ and the remaining vertices in $H_{[\ell]}$ are partitioned into $\mathcal{X}^*$.
By Claim~\ref{cl:rearrangement_of_H}~(c),(d) and \eqref{eq: size n^*}, the set $X^*_i$ in $\mathcal{X}^*$ has size very close to $n^*_i=|V^*_i|$, but they are not exactly same.
Our plan is that for each $i\in [r]\setminus[k]$, we move $n^*_i-|X^*_i|$ vertices from $V^*_i$ to $V^*_j$ for some $j\in [k]$ so that we have exactly the desired number of vertices in each $V^*_i$ and $V^*_j$.

Since $H$ and $G$ have exactly the same number of vertices, Claim~\ref{cl:rearrangement_of_H}~(c),(d) and \eqref{eq: size n^*} imply that each of the following summands is nonnegative and we have the equality:
$$\sum_{j\in [k]} |X^*_j| - n^*_j = \sum_{i\in [r]\setminus [k]} n^*_i - |X^*_i|.$$
Hence, we consider an auxiliary bipartite multi-graph $Aux$ with bipartition $([k],[r]\setminus[k])$ such that the degree of vertex $a\in [k], b\in [r]\setminus [k]$ is $d_{Aux}(a)= |X^*_a|-n^*_a$ and $d_{Aux}(b) = n^*_b-|X^*_b|$. Note that $|E(Aux)|\leq \eta^{1/5} n$ by Claim~\ref{cl:rearrangement_of_H}~(c),(d) and \eqref{eq: size n^*}.

For each edge $ij\in Aux$ with $j\in [k]$, we consider the path $P(i,j)$ in $D$ given in \eqref{eq: dipath in D}. 
Assume that $P(i,j)= p_1\dots p_\ell$. Consider an edge $p_{a}p_{a+1}$ in $P(i,j)$.
For each $j'\in N_{Q}(p_{a+1})$, the fact $p_{a}p_{a+1}\in E(D)$ implies that $(V_{p_{a}},V_{j'})$ is an $(\eps',d)$-regular pair in $G$. Hence, by Lemma~\ref{lem:large_super_regular}, at least $(1- k \eps'^{1/2} ) n^*_{p_{a}}$ vertices in $V^*_{p_a}$ has at least $(d-\eps'^{1/2})|V^*_{j'}|$ neighbors in $V^*_{j'}$ for all $j'\in N_Q(P_{a+1})$. We choose one such vertex and move it to $V^*_{p_{a+1}}$ then the super-regularity of $G$ with respect to the clique $Q$ is not ruined much.
We repeat this for all edges in $P(i,j)$ for all $ij\in E(Aux)$ while updating the set $V^*_{i}$.
As each path $P(i,j)$ has length at most $r$, we have to do this for at most $r|E(Aux)|\leq r\eta^{1/5}n\leq \eta^{1/6}n < \eps'\frac{n}{r}$ times, the regularity of $G$ between each regular pair does not get ruined, so we can complete this task until the end. Moreover, by the definition of $Aux$, at the end of this process,  we obtain the new partition $\tilde{V}_i$ of $V(G)\setminus \phi'_{\ell}(V(H^*))$ with the property that $|X^*_i| = |\tilde{V}_i|$. By \Cref{lem:tiny_part_good_enough,lem:adding_tiny_parts}, we have the following.
\begin{equation}
    \begin{minipage}[c]{0.9\textwidth}
    \begin{enumerate}
    \item[$\bullet$] For each $ij\in E(Q)$, the pair $(\tilde{V}_i,\tilde{V}_j)$ is $(\eps'^{1/4},d/2)$-super-regular. \label{after_move_super_regular}
    \item[$\bullet$] For each $ij\in E(R)$, the pair $(\tilde{V}_i,\tilde{V}_j)$ is $(\eps'^{1/4},d/2+)$-regular.
    \item[$\bullet$] For each $i\in [r]$, we have $|\tilde{V}_i|=|X^*_i| = (1\pm \eps'^{1/10}) n_i$ 
\end{enumerate}
    \end{minipage}
\end{equation}
Now, for each $i\in [r']$, we consider the graph $H_{[\ell]}[X'_{Q^i}]$ and $G[\tilde{V}_{Q^i}]$ where $X'_{Q^i}=\bigcup_{q\in Q^i} X'_q$ and $\tilde{V}_{Q^i} = \bigcup_{q\in Q^i} \tilde{V}_q$.
Let $T$ be the union of $T_{j}$ and $T'_{j+1}$ where $H_j$ is assigned to $Q^i$ in the application of Claim~\ref{cl:rearrangement_of_H}, i.e. $s_j=i$.
Recall that for each vertex $x\in T$, we defined the target set $A_x$ in \eqref{available_vertex_C}.
 
If $i=r'$, then the graph $H_\ell$ is in $X'_{Q^i}$. Then  we define the hypergraph $\M$ to be the hypergraph $\M$ in the assumption of \Cref{lem:technical_bandwidth} (otherwise, define it to be empty hypergraph). Apply \Cref{lem:technical_blow} with the following sets: 
\begin{center}
\begin{tabular}{ |c|c|c|c|c|c|c|c| }
\hline
 sets/graphs & $G[\tilde{V}_{Q^i}]$ & $H_{[\ell]}[X'_{Q^i}]$ & $\M$ & $\tV_i$ & $U_i$ & $T$ & $A_x$ \\ 
 \hline
 playing the role of & $G$ & $H$ & $\M$ & $V_i$ & $U_i$ & $T$ & $A_x$ \\ \hline
\end{tabular}
\end{center}
It can be easily checked that the hypotheses of Lemma~\ref{lem:technical_blow} hold with the following parameters in all of the above cases. 
\begin{center}
\begin{tabular}{ |c|c|c|c|c|c|c|c|c|c| }
\hline
 parameters/numbers & $2n/r$ & $2r\a/\h$ & $\e'^{1/4}$ & $(d/4)^\D$ & $1/4$ & $d$ & $\D$ & $k$ & $s$ \\ 
 \hline
 playing the role of & $n$ & $\g$ & $\e$ & $\a$ & $\b$ & $d$ & $\D$ & $k$ & $s$ \\ \hline
\end{tabular}
\end{center}
The conditions \ref{blowup_assumption_0} and \ref{blowup_assumption_1} follow from \eqref{after_move_super_regular} using the identity function $g$. \ref{blowup_assumption_2} and \ref{blowup_assumption_3} trivially hold. The condition \ref{blowup_assumption_4} follows from \eqref{available_vertex_C}.

The conclusion of \ref{blowup_conclusion_1} ensures that the vertices are embedded in the corresponding available sets only (in other words, the adjacency relations are preserved in the embedding). 
This yields an embedding of $H_{[\ell]}[X'_{Q^i}]$ into $G[\tilde{V}_{Q^i}]$ satisfying \ref{blowup_conclusion_1} and \ref{blowup_conclusion_2}.
Repeating this for all $i\in [r']$ together with $\phi'_{\ell}$ yields a desired embedding of $H_{[\ell]}$. As each embedding of $H_{[\ell]}[X'_{Q^i}]$ satisfies the property \ref{blowup_conclusion_1}, they together with $\phi'_{\ell}$ provide an embedding of $H$ into $G$, and the property \ref{blowup_conclusion_2} and Claim~\ref{cl:rearrangement_of_H}-(e) ensure \ref{bandwidth_conclusion_2}.
This finishes the proof of \Cref{lem:technical_bandwidth}.

\end{proof}

\section{Proof of Theorem~\ref{main}}
In this section, we prove Theorem~\ref{main}. 
We will follow the steps described in Section~\ref{sec: proof sketch}. 
Let $V$ be the common vertex set of the graphs in $\mathcal{G}$. 

We first claim that we may assume that $H$ has no isolated vertices.
If $H$ has $n'$ non-isolated vertices, then $n'\geq h/\D$ is sufficiently large. Then we take a subset $V'\subseteq V$ of size exactly $n'$ where $\delta(G_i[V']) \geq (1-\frac{1}{k}+\frac{1}{2}\eps) n'$ for all $i\in [h]$. This will allow us to be able to deal with the cases when $H$ has no isolated vertices. 
In order to show that such a choice of $V'$ exists, we consider a random subset $V'\subseteq V$ of size exactly $n'$ and analyze the probability of getting $\delta(G_i[V']) \geq (1-\frac{1}{k}+\frac{1}{2}\eps) n'$ for all $i\in [h]$. 
We call a vertex $v\in V$ \emph{bad} if it is chosen to be in $V'$ and less than $(1-\frac{1}{k}+\eps/2)n'$ neighbors of it in $G_i$ are chosen for some $i\in [h]$. 
We claim that a vertex $v$ is bad with the probability at most $\frac{n'}{n}\cdot h \cdot e^{-\eps^2n'/100}$. 
Indeed, a fixed vertex $v$ is chosen with probability $n'/n$. Conditioning on that, we consider the hypergeometric distribution of choosing $n'-1$ vertices from $n-1$ vertices. Applying Lemma~\ref{lem:concentration} yields us the desired bound.
Consequently, the linearity of expectation shows that the expected number of bad vertices is at most $n'\cdot h \cdot e^{-\eps^2n'/100} <1$ as $h/\Delta\leq n'$ and $1/h\ll \eps, 1/\Delta$. This shows that there is a choice of $V'$ with no bad vertices, yielding the desired choice of $V'$. Therefore, we may assume that $H$ has no isolated vertices, which will be crucial for Step~2-2.

We start with partitioning $H$ into $H_0,\dots, H_{m+1}$ using an $\alpha n$-bandwidth partition.
We choose further constants $\gamma,\beta_0,\beta_1,\beta_2,r'_0,r'_1, \eps_0,d',d$ so that the followings hold: $$0<\frac{1}{h_0} \ll \alpha \ll \gamma\ll \beta_2\ll \beta_1 \ll \beta_0 \ll  \frac{1}{r'_1} \ll \frac{1}{r'_0} \ll \eps_0 \ll d'\ll d \ll \eps, \frac{1}{k}, \frac{1}{\Delta} < 1.$$
As $h\geq h_0$  and $n> h/\Delta$, we have $1/n \ll \alpha$.
By considering a $k$-coloring $c$ of $H$, we consider an $\alpha n$-bandwidth partition $W_1,\dots, W_t$ of $H$ with respect to the coloring. We have that $k$ divides $t$ and some sets $W_i$ may be empty. 
In fact, some $W_i$ might be empty, but at least one set among $k$ consecutive sets $W_{i},\dots, W_{i+k-1}$ are non-empty from the definition of the bandwidth partition (see Lemma~\ref{lem: bandwidth partition}). Let $Y_1,\dots, Y_k$ be the color classes of the coloring $c$. In particular, $W_i\subseteq Y_j$ if $i\equiv j~({\rm mod}~k)$ for some $j\in [k]$.

In order to partition $H$, we choose indices $t_0,t_1,\dots, t_m, t_{m+1}=t$ for an appropriate $m$ so that the following holds.
\begin{equation}\label{eq: choice of Hsss}
    \begin{minipage}[c]{0.9\textwidth}
\begin{enumerate}
    \item[(a)] $|E(H[W_{[t_0]}])| = (\beta_0\pm \alpha^{1/2}) n$ and $|E(H[W_{[t_{m}+1,t]}])| = (\beta_2 \pm 2 k\Delta \gamma) n$.
    \item[(b)] $|E(H[W_{[t_{i}+1,t_{i+1}]}])| = (\gamma \pm \alpha^{1/2})n$ for each $i\in [m]$.
    \item[(c)] $k$ divides each of $t_0,\dots, t_m$.
\end{enumerate}
\end{minipage}
\end{equation}
This choice is possible as $\alpha \ll \gamma, \beta_0,\beta_2 \ll 1/\Delta$ and $\Delta(H)\leq \Delta$. Furthermore, it ensures that $\gamma^{-1}/4\leq m\leq 4\Delta \gamma^{-1}$. Let 
$$H_0 = H[W_{[t_0]}],\enspace H_i= H[W_{[t_{i-1}+1,t_i]}] \enspace \text{for each $i\in [m+1]$}.$$ 
For each $s\in [m+1]\cup \{0\}$, let $n_i = |V(H_i)|$ and let $H_{[s]}$ denote the graph $H[W_{[t_s]}] = H[\bigcup_{0\leq i\leq s} V(H_i) ]$. 
Since $H$ has no isolated vertices, we have $n_i\geq \gamma n/\Delta$. \newline

\noindent{\bf Step 1. Find a color absorber {\boldmath $(\varphi_0(H_0),A,B,C)$}.}
In this step, we will find an embedding $\varphi_0$ of $H_0$ into $\mathcal{G}^{\eps/2}$.
This ensures that each edge in $\varphi_0(H_0)$ has many possible colors. This together with Lemma~\ref{absorption} will yield a desired color absorber.
Let
$$G^0 = \mathcal{G}^{\eps/2}.$$
Then by Lemma~\ref{lem_mindeg}, we have
$$\delta(G^0) \geq \left(1-\frac{1}{k}+\frac{\eps}{2}\right)n.$$
We will use Lemma~\ref{lem:technical_bandwidth} to find an embedding $\varphi_0 : H_0\rightarrow G^0$.
For this, we will first find a regularity partition $(R^0, \mathcal{V}^0)$ for $G^0$ and find an appropriate subset $V^0\subseteq V$ and embed $H_0$ into $G^0[V^0]$.

We apply the regularity lemma (Lemma~\ref{regularity}) to $G^0$ with the parameters $\eps_0, r'_0, r'_1,d$ playing the roles of $\eps, t, M,d$, respectively to obtain an $\eps_0$-regular equitable partition $\mathcal{V}=(V_1, \cdots, V_{r_0})$ with $r'_0\leq r_0\leq r'_1$ with the reduced graph $R^0$ on the vertex set $[r_0]$ and $\delta(R^0)\geq \left(1-\frac{1}{k} +\frac{\eps}{4}\right)r_0$ and $k$ divides $r_0$. 
Since the minimum degree condition on $R^0$ ensures the existence of many $k$-cliques in $R^0$, we can assume that $[k]$ is a $k$-clique in $R^0$.

We choose a set $V^0\subseteq V$ with $|V^0|= n_0$ satisfying the following where we let $V_i^0 = V_i \cap V^0$, $U^0= V\setminus V^0$, $U_i^0 = V_i\setminus V_i^0$, $\mathcal{U}^0 = (U^0_1,\dots, U^0_r)$ and $\mathcal{V}^0 = (V^0_1,\dots, V^0_{r_0})$.

\begin{equation}\label{eq: U0 property}
    \begin{minipage}[c]{0.9\textwidth}
\begin{enumerate}[leftmargin=*]
    \item[(a)] $|V^0_i|=(1 \pm \eps_0^{1/4}) \frac{n_0}{r_0}$ for all $i \in [r_0]$. 
    \item[(b)] $G^0[V^0]$ admits an $(\eps_0^{1/10},d)$-regular partition $(R^0,\mathcal{V}^0)$, and $[k]$ forms a clique in $R^0$.
    \item[(c)] For all $ij\in E(R^0)$, the pairs $(V^0_i,U^0_j)$ and $(U^0_i,U^0_j)$ are $(\eps_0^{1/10}, d+)$-regular pair in~$G^0$.
    \item[(d)] For each $i\in [h]$, we have $\delta(G_i[V^0]) \geq \left(1-\frac{1}{k} + \frac{\eps}{2}\right) |V^0|$ and $\delta(G_i[U^0])\geq \left(1-\frac{1}{k} + \frac{\eps}{2}\right)|U^0|$.
\end{enumerate}
    \end{minipage}
\end{equation}
Indeed, to show that such a set exists, we choose a random set $V^0$ of size exactly $n_0$ from $V$. Then using Lemma~\ref{lem:randomset_reg}, (a), (b) and (c) holds with probability at least $1-o(1/n)$. Using Lemma~\ref{lem:concentration}, we can show that (d) holds with probability at least $1-o(1/n)$ by considering random variables measuring the degree on the subgraph $G_i[V^0]$ for each vertex $v\in V$ and $i\in [h]$ and applying the union bound.  
 
Fix one such choice of $V^0$ as above. We need to ensure that we can embed $H_1$ into $U^0$ after embedding $H_0$ into $V^0$. For this purpose, we need to make sure that the last few vertices in $H_0$ are embedded in such a way that the connection between $H_0$ and $H_1$ can be established. 
As $[k]$ forms a clique in $R^0$, we plan to embed the last $k$ sets $W_{t_0-k+1},\dots, W_{t_0}$ into $V^0_1,\dots, V^0_k$.
Then, for each $x\in W_{[t_0+1,t_{0}+k-1]}$, the neighborhood $N_x = N_{H}(x)\cap W_{[t_0]}$ must be embedded in such a way that the vertices in $\varphi_0(N_x)$ have a large common neighborhood in an appropriate set $U^0_\ell$ within the graph $G^0$. This can be ensured exploiting \ref{bandwidth_conclusion_2} of Lemma~\ref{lem:technical_bandwidth}.
For this, we let
$$\M^0 = \{ N_x: x\in W_{[t_0+1,t_{0}+k-1]} \}$$ 
be a multi-hypergraph on the vertex set $W_{[t_0-k+1,t_0]}$ of rank at most $\Delta$ and $\Delta(\M^0)\leq \Delta$. For each $x\in W_{t_0+\ell}$ with $\ell\in [k-1]$, we let $f^0(N_x) = \ell$.
Now we apply Lemma~\ref{lem:technical_bandwidth} to embed $H_0$ into $G^0[V^0]$ with the following parameters.
\begin{center}
\begin{tabular}{ |c|c|c|c|c|c|c|c|c| }
\hline
 parameters/objects & $n_0$ & $\eps/4$ & $d$ & $d'$ & $\eps_0^{1/10}$ & $r_0$ & $1/r_0^2$ & $k-1$\\ 
 \hline
 playing the role of & $n$ & $\eps$ & $d$ & $d'$ & $\eps'$ & $r$ & $\eta$ & $s$\\ \hline \hline
  parameters/objects &  $\alpha \frac{n}{n_0}$ & $[k]$ & $U_i^0$ & ${U_i}^0$ & $\varnothing$ & $\{W_1,\dots, W_{t_0}\}$ & $G^0[V^0\cup U^0_{[k-1]}]$ & $\M^0$\\ 
 \hline
 playing the role of &  $\alpha$ & $K$ & $U_i$ & $V_i$ & $\mathcal{I}$ & $\mathcal{W}$ & $G'$ & $\M$\\
 \hline
\end{tabular}
\end{center}
Then, as $\frac{n}{n_0}$ is between $\frac{1}{2}(\beta_0 \pm k \alpha)^{-1}$ and $\Delta(\beta_0\pm k\alpha)^{-1}$, the assumption $\alpha \ll \beta_0$ implies that $\alpha \frac{n}{n_0}$ is small enough. 
So, the hierarchy for the application of Lemma~\ref{lem:technical_bandwidth} is met.
\ref{cond: G regularity} holds by the above condition \eqref{eq: U0 property}-(b) and \ref{cond: M hypergraph} holds by the definition of $\M$ and the fact that $\Delta(H)\leq \Delta$. 
\ref{cond: M regularity} holds by \eqref{eq: U0 property}-(c) as $[k]$ forms a $k$-clique in $R^0$. 
The last condition \ref{cond: target sets} is vacuous as we do not use any fragmented interval in this application. Thus Lemma~\ref{lem:technical_bandwidth} yields an embedding $\varphi_0:H_0 \rightarrow G^0[V^0]$ such that for every $\ell\in [k-1]$ and $x \in W_{t_0+\ell}$, we have 
\begin{align}\label{eq: matching embedding 0}
    |A_x| \geq d' \left|U_{\ell}^0\right|,
\end{align}
where $A_x = \bigcap_{y\in N_H(x)\cap W_{[t_0]}} N_{G^0}(\varphi(y))\cap U_{\ell}^0$.
Finally, we apply Lemma~\ref{absorption} with $n,h,\eps/2, e(H_0)/n$ and $\beta_1$ playing the roles of $n,h,\eta,\beta$ and $\gamma$, respectively.
Here, $e(H_0)\leq \Delta n_0$ and $H$ has no isolated vertices, so $e(H_0)/n$ is between $\beta_0/4$ and $2\Delta \beta_0$, thus the application of Lemma~\ref{absorption} is possible.
This yields two disjoint sets $A, C\subseteq [h]$ satisfying the following.
\begin{equation}\label{eq: color absorption}
    \begin{minipage}[c]{0.9\textwidth}
\begin{enumerate}[leftmargin=*]
    \item[(a)] $e(H_0) -|A| = \beta_1 n$ and $|C| = (\beta_0\pm \alpha^{1/2}) n$.
    \item[(b)] For every subset $C' \subseteq C$ of size exactly $e(H_0) -|A|=\beta_1 n$, there exists a bijection ${\lambda_0:E(\varphi_0(H_0))\rightarrow [A\cup C']}$, yielding an $\mathcal{G}_{A\cup C'}$-transversal isomorphic to $H_0$. 
\end{enumerate}
\end{minipage}
\end{equation}
Let $B= [h]\setminus (A\cup C)$. The tuple $(\varphi_0(H_0),A,B,C)$ is a desired color absorber. \newline

\noindent{\bf Step 2. Embed {\boldmath $H_1,\dots, H_{m+1}$} using all colors in {\boldmath $B$}.}
We now embed $H_1,\dots, H_m, H_{m+1}$ into $\binom{V}{2}$. The color absorber $\varphi_0(H_0)$ will stay uncolored until the end of the proof, but we will actually assign colors for the images of each of $H_1,\dots, H_m, H_{m+1}$. While doing this, we will only use the colors in $B\cup C$. 
All embedding procedures are almost the same, except that embedding $H_{m+1}$ is a bit different as we want to exhaust all colors in $B$ while embedding $H_{m+1}$.

When we embed $H_s$, we will first choose a set $B_s$ of colors and embed $H_s$ into the graph $\mathcal{G}_{B_s}^{\eps/8}$. While doing that, we use Lemma~\ref{lem: connection} to build connection between $H_{s-1}$ and $H_s$ using the colors in $B_{s-1}$ and $B_s$. For this purpose, we let $B_0 = B$.

Assume that for $s\in [m]\cup \{0\}$, we have an embedding $\varphi_s$ of $H_{[s]}$ into $\binom{V}{2}$ and a graph $G^s$ and a set $U^s$ and a partial coloring $\lambda_s$ of $\varphi_s(H_{[s]})$ satisfying the following where $U^{s} = V\setminus \varphi_s(H_{[s]})$ and $B_s$ is a subset of $B\cup C$ and $G^{s} = \mathcal{G}_{B_s}^{\eps/8}$. 
\stepcounter{propcounter}
\begin{enumerate}[label = {\bfseries \emph{\Alph{propcounter}\arabic{enumi}}}]
\item The partial coloring $\lambda_s$ is an injection from $E(H_{[s]})\setminus E(H_0)$ to $B\cup C$ such that for every $e\in E(H_{[s]})\setminus E(H_0)$, we have $\varphi_s(e) \in G_{\lambda_s(e)}$. Let $\Lambda_s$ be the image of $\lambda_s$. In addition, if a color in $B_s$ is used, then it is used in coloring an edge of $E( V(H_{s-1})\cup V(H_{s}))$.\label{eq: cond1 for step 2}
\item $|B_s|\geq \gamma^{1/2}n$  \label{eq: cond2 for step 2}
\item $\Lambda_s \subseteq B$ if $|B\setminus \Lambda_{s-1}| > \gamma^{1/2}n$. \label{eq: cond3 for step 2}
\item For each $i\in [h]$, we have $\delta(G_i[U^s])\geq \left(1-\frac{1}{k}+ \frac{\eps}{2} - s\gamma^2\right)|U^s|$. \label{eq: mindeg GiU}
\item $G^s[U^s]$ admits an $(\eps_0^{1/20}, d/2)$-regular partition $(R^s,\mathcal{U}^s)$ with $\mathcal{U}^s = (U^s_1,\dots, U^s_{r_s})$ and ${r'_0\leq r_s\leq r'_1}$, where $r_s = |V(R^s)|$. In addition, $[k]$ forms a clique in $R^s$.\label{eq: Us regularity}
\item For each $\ell\in [k-1]$ and $x\in W_{t_s+\ell}$, we have $|A_x|\geq  d' |U^s_\ell|$ where $$A_x = \bigcap_{y\in N_H(x)\cap W_{[t_s]}} \left(N_{G^s}(\varphi_s(y))\cap U^s_{\ell}\right).$$\label{eq: target sets in step 2}
\end{enumerate}
Note that \eqref{eq: U0 property} implies that $\varphi_0$ satisfies the above six properties.
We now show that assuming that $\varphi_s$ satisfies the above properties, we can find an embedding $\varphi_{s+1}$ of $H_{[s+1]}$ satisfying the above six properties. This ensures that we can repeat this until we obtain $\varphi_{m+1}$ embedding the entire graph. As mentioned before, the case $s=m$ will be dealt with extra care as we plan to exhaust all colors in $B$ in this case completely. 

If $|B\setminus \Lambda_s|> \gamma^{1/2}n$, then choose a set $B_{s+1}\subseteq B\setminus \Lambda_s$ of size exactly $\gamma^{1/2} n$; thus, \ref{eq: cond2 for step 2} trivially holds with $B_{s+1}$. Otherwise, we choose $B_{s+1}$ to be $(B\cup C)\setminus \Lambda_s$. Since $s\leq m$, \eqref{eq: choice of Hsss}-(a) and \eqref{eq: color absorption}-(a) ensure that $$|B_{s+1}| \geq e(H_{m+1})+ e(H_0)-|A|\geq \beta_1 n\geq \gamma^{1/2} n,$$ which, in turn, ensures that \ref{eq: cond2 for step 2} holds with $B_{s+1}$.

Let $B'= B\setminus \Lambda_s$ be the set of remaining unused colors in $B$.
If $s=m$, then we additionally have the following properties.
\begin{align} \label{eq: property B'}
    B'\subseteq B_{m+1}, \enspace |B'|\leq \gamma^{1/2}n  \enspace \text{ and } \enspace 
    |B_{m+1}| \geq \beta_1 n.
\end{align}
Indeed, if 
$s=m$, then the used colors $\lambda_s(\varphi_{s}(H_s))$ contains all colors in $[h]$ except at most $$ |E(H_0)| + |E(H_{[t_{m}+1, t]} )| + |E(H[W_{[t_{m}]},W_{[t_{m}+1,t_{m}+k]} ])|$$
colors. By \eqref{eq: choice of Hsss}-(a) and \eqref{eq: color absorption}-(a), this number is at most $|A| + \beta_1 n + \beta_2 n + 6k\Delta \gamma n \leq |A|+|C| + 10 k\Delta \gamma n$. That means that we have $|B\setminus \Lambda_{s}| \leq 10 k\Delta \gamma n  < \gamma^{1/2}n$, hence $B_{s+1}$ is chosen as $(B\cup C)\setminus \Lambda_s$, which contains at least $\beta_1 n$ colors.
If $s=m$, then we have an additional goal of exhausting all the colors in $B'$ while embedding $H_{m+1}$.

We will embed $H_{s+1}$ into the following graph.
$$G^{s+1} = \mathcal{G}^{\eps/8}_{B_{s+1}}.$$ Then Lemma~\ref{lem_mindeg} together with \ref{eq: mindeg GiU} implies
\begin{align*}
    \delta(G^{s+1}[U^s])\geq \left(1-\frac{1}{k}+ \frac{\eps}{4}\right)|U^s|.
\end{align*}
In order to use Lemma~\ref{lem: connection} to embed $H_{[t_s+1,t_{s+1}]}$ into $G^{s}\cup G^{s+1}$, we need regularity partitions of $G^s$ and $G^{s+1}$, which we obtain next. 

Apply Lemma~\ref{regularity} to $G^{s+1}[U^s]$ to obtain a regularity partition $(R^{s+1},\mathcal{V})$ with $\mathcal{V}=(V_1,\dots, V_{r_{s+1}})$ satisfying the following.
\begin{equation}\label{eq: middle V}
    \begin{minipage}[c]{0.9\textwidth}
\begin{enumerate}[leftmargin=*]
    \item[$\bullet$] $|V_i|= \frac{1}{r_{s+1}}|U^s| \pm 1$ for all $i \in [r_{s+1}]$ and $k$ divides $r_{s+1}$.
    \item[$\bullet$] $G^{s+1}[U^s]$ admits an $(\eps_0,d)$-regular partition $(R^{s+1},\mathcal{V})$ with $|R^{s+1}|=r_{s+1}$ and ${r'_0\leq r_{s+1}\leq r'_1}$.
    \item[$\bullet$] $\delta(R^{s+1}) \geq \left(1- \frac{1}{k} +\frac{\eps}{8}\right) r_{s+1}$. Moreover, we assume $[k]$ forms a clique in $R^{s+1}$.
\end{enumerate}
    \end{minipage}
\end{equation}
Similarly as we did in \eqref{eq: U0 property}, by taking a random subset of $V^{s+1}$ of size exactly $n_{s+1}= |V(H_{s+1})|$, we can show that there exists a subset $V^{s+1}\subseteq U^s$ satisfying the following where $U^{s+1}= U^s\setminus V^{s+1}$, $V_i^{s+1} = V_i \cap V^{s+1}$, $U_i^{s+1} = V_i \setminus V^{s+1}$ and $\mathcal{V}^{s+1} = (V^{s+1}_1,\dots, V^{s+1}_{r_{s+1}})$.
\begin{equation}\label{eq: Us+1 property}
    \begin{minipage}[c]{0.9\textwidth}
\begin{enumerate}[leftmargin=*]
    \item[(a)] $|V^{s+1}_i|=(1 \pm \eps_0^{1/20}) \frac{ n_{s+1}}{r_{s+1}}$ for all $i \in [r_{s+1}]$ and $k$ divides $r_{s+1}$.
    \item[(b)] $G^{s+1}[V^{s+1}]$ admits an $(\eps_0^{1/20},d/2)$-regular partition $(R^{s+1},\mathcal{V}^{s+1})$.
    \item[(c)] For $ij\in E(R^{s+1})$, the pair $(V^{s+1}_i,U^{s+1}_j)$ and $(U^{s+1}_i,U^{s+1}_j)$ are $(\eps_0^{1/20}, d/2+)$-regular pairs in~$G^{s+1}$.
    \item[(d)] For each $i\in [h]$, we have $\delta(G_i[V^{s+1}]) \geq \left(1-\frac{1}{k} + \frac{\eps}{2} - (s+1)\gamma^2\right) |V^{s+1}|$ and $\delta(G_i[U^{s+1} ])\geq \left(1-\frac{1}{k} + \frac{\eps}{2} - (s+1)\gamma^2 \right)|U^{s+1}|$.
    \item[(e)] The graph $G^s$ admits an $(\eps_0^{1/150},d/2)$-regular partition $(R^s, \widetilde{\mathcal{U}}^s)$ with $\widetilde{\mathcal{U}}^s=(\widetilde{U}^s_1,\dots,\widetilde{U}^s_{r_s})$ with $\widetilde{U}^s_{i} = U^s_{i} \cap V^{s+1}_i$.
\end{enumerate}
    \end{minipage}
\end{equation}
The last property (e) can also be shown from \ref{eq: Us regularity} using Lemma~\ref{lem:randomset_reg}.
If $s=m$, then we have $V^{s+1}= U^s$ and $U^{s+1}=\varnothing$, thus the above properties become immediate from \eqref{eq: middle V}, \ref{eq: mindeg GiU} and  \ref{eq: Us regularity}. \newline

Since we wish to use Lemma~\ref{lem: connection} and Lemma~\ref{lem:technical_bandwidth}, we divide $H^{s+1}$ into the following two graphs. 
$$H^{\rm con} = H[W_{[t_{s}+1,t_{s}+4k^3]}] \enspace \text{and} \enspace H^{\rm band}= H[W_{[t_s+4k^3+1,t_{s+1}]}].$$

\noindent {\bf Step 2-1. Embedding {\boldmath $H^{\rm con}$} by Lemma~\ref{lem: connection}.}
Similarly to the definition of $\M^{0}$ before, we want to again define $\M^{s+1,1}$ to establish connection between $H^{\rm con}$ and $H^{\rm band}$
so that we can later embed $H^{\rm band}$.
Since $[k]$ forms a clique in $R^{s+1}$, we plan to embed the last $k$ sets $W_{t_s+4k^3-k+1},\dots, W_{t_s+4k^3}$ into $V^{s+1}_1,\dots, V^{s+1}_k$.
Then, for each $x\in W_{[t_s+4k^3+1,t_s+4k^3+k-1]}$, consider $N_x = \{N_{H}(x)\cap W_{[t_s+4k^3]}\}$.
We let
$$\M^{s+1,1} = \{ N_x: x\in W_{[t_s+4k^3+1,t_s+4k^3+k-1]} \}$$ 
be a multi-hypergraph on the vertex set $W_{[t_s+4k^3-k+1,t_s+4k^3]}$ of rank at most $\Delta$ and $\Delta(\M^0)\leq \Delta$. For each $x\in W_{t_s+4k^3+\ell}$ with $\ell\in [k-1]$, we let $f(N_x) = \ell$.
Now, we apply Lemma~\ref{lem: connection} to embed $H^{\rm con}$ to the two graphs $G^{s}[V^{s+1}]$ and $G^{s+1}[V^{s+1}]$ with the following parameters and objects
\begin{center}
\begin{tabular}{ |c|c|c|c|c|c|c|c|c|c|c|c|c| }
\hline
 parameters/objects & $G^{s}[V^{s+1}]$& $\widetilde{\mathcal{U}}^s$ &$G^{s+1}[V^{s+1}]$& $\mathcal{V}^{s+1}$ &$n_{s+1}$ & $\eps/8$ & $d/2$ & $\alpha \frac{n}{n_{s+1}}$\\ 
 \hline
 playing the role of & $G^1$ & $\mathcal{V}^1$ &$G^2$ & $\mathcal{V}^2$ & $n$ & $\eps$ & $d$ & $\alpha$\\ \hline \hline
 parameters/objects & $d'$ & $\eps_0^{1/150}$ &  $r'_1$ & $A_x$ & $\M^{s+1,1}$ & $f$ & $4k^3$ & $W_{t_s+i}$\\ 
 \hline
 playing the role of  & $d', \beta$ & $\eps_0$ & $r$ & $A_x$ & $\M$ & $f$ & $l$ & $W_i$\\
 \hline
\end{tabular}
\end{center}
Indeed, this application is possible. Indeed, it is easy to check that the hierarchy of the constants is satisfied and
\eqref{eq: Us+1 property}-(a),(b),(d),(e) imply that all conditions for Lemma~\ref{lem: connection} are met.
This provides an embedding $\varphi'$ of $H^{\rm con}$ into $G^{s}\cup G^{s+1}$ so that the following two holds.
\begin{enumerate}[leftmargin=*]
    \item[(a)] For each $x\in W_{t_{s}+j}$ with $j\in [k-1]$, $\varphi'(x)\in A_x$, where the set $A_x$ is defined in \ref{eq: target sets in step 2}.
    \item[(b)] For each $j\in [k-1]$ and $x\in W_{t_{s}+4k^3+j}$, we have the following where ${A_x = \bigcap_{y\in N_H(x)} N_{G^{s+1}}(\varphi'(y)) \cap V^{s+1}_{j}}$.
    \begin{equation}\label{eq: size of Ax}
    |A_x| \geq d' |V^{s+1}_j|.
    \end{equation}
\end{enumerate}
The property $(b)$ comes from the definition of $\M^{s+1,1}$.
Let $\varphi'_{s} = \varphi_s \cup \varphi'$, then (a) ensures that $\varphi'_s$ maps all edges of $H$ between $W_{[t_s]}$ and $W_{[t_s+1,t_s+ k-1]}$ into $G^s$. In this embedding, the edges in $\varphi'_s(H[W_{[t_s+4k^3]}])\setminus \varphi_s(H_{[s]})$ are not yet colored.

For each edge $e\in (\varphi'_s(H[W_{[t_s+4k^3]}])\setminus \varphi_s(H_{[s]}))\cap E(G^{s})$, we choose a color $\lambda'_{s+1}(e)$.
Since $G^{s} = \mathcal{G}_{B_s}^{\eps/8}$, there are at least $(\eps/8)|B_s|$ colors $i\in B_s$ such that $e\in G_i$. The property \ref{eq: cond1 for step 2} ensures that at most $e(V(H_{s-1}\cup H_{s}))\leq \Delta |V(H_s)|+\Delta |V(H_{s-1})|\leq 4\Delta \gamma n$ colors in $B_s$ are already used to color an edge in $\varphi_{s}(H_{[s]})$.
Hence, there are still $(\eps/8)|B_s| - 2\Delta \gamma n> \gamma^{3/4}n > e(H^{\rm con})$ colors are available for each $e\in \varphi'(H^{\rm con} \cap G^{s})$. 
We greedily color those edges to obtain a partial rainbow coloring of the embedded graph.
Again, in the same way, we color each edge $e\in \varphi'(H^{\rm con}) \cap G^{s+1}$ by a color in $B_{s+1}$ in a greedy way. This uses at most $e(H^{\rm con})\leq 4k^3\Delta \alpha n$ colors from $B_{s+1}$, making most of the colors in $B_{s+1}$ still available. This together with the previous partial coloring $\lambda_s$ yields a new partial coloring $\lambda'_{s}$ coloring all edges in $\varphi'_{s}(H[W_{[t_s+4k^3]}])\setminus \varphi'_{s}(E(H_0))$. 

Let $B'' = B'\setminus \lambda'_{s}(H[W_{[t_s+4k^3]}])$ be the set of remaining available colors in $B'$.
Let $V'' = \varphi'(V(H^{\rm con}))$, then $|V''|\leq 4k^3 \alpha n$. 
Since $\alpha \ll 1/r_{s+1}$ and $1/r_{s+1}\ll \eps_0$, after deleting all vertices in $V''$ from each of $V^{s+1}_i$, the properties $(a),(b),(c),(e)$ in \eqref{eq: Us+1 property} hold when replacing $\eps_0$ with $2\eps_0$.
The property \eqref{eq: Us+1 property}-(c) holds with 
\begin{align}\label{eq: min deg for next}
    \delta(G_i[V^{s+1}\setminus V'']) \geq \left(1- \frac{1}{k} + \frac{\eps}{2} - (s+1)\gamma^2 \right)|V^{s+1}| - 4k^2\alpha n\geq \left(1- \frac{1}{k} + \frac{\eps}{4} \right)|V^{s+1}|.
\end{align}
We now embed $H^{\rm band}$ into $V^{s+1}\setminus V''$. 
We plan to embed the last $k$ sets $W_{t_{s+1}-k+1},\dots, W_{t_{s+1}}$ into $V^{s+1}_1\setminus V'',\dots, V^{s+1}_k\setminus V''$ as we assumed $[k]$ forms a clique in $R^{s+1}$.
Then, again if $s<m$, then for each $x\in W_{[t_{s+1}+1,t_{s+1}+k-1]}$, consider the neighborhood $N_x = \{N_{H}(x)\cap W_{[t_s]}\}$ and we let
$$\M^{s+1,2} = \{ N_x: x\in W_{[t_{s+1}+1,t_{s+1}+k-1]} \}$$ 
be a multi-hypergraph on the vertex set $W_{[t_{s+1}-k+1,t_{s+1}]}$ of rank at most $\Delta$ and $\Delta(\M^{s+1,2})\leq \Delta$. For each $x\in W_{t_{s+1}+\ell}$ with $\ell\in [k-1]$, we let $f(N_x) = \ell$.
If $s=m$, then $\M^{s+1,2}$ is an empty hypergraph.\newline

\noindent {\bf Step 2-2. Assign target sets {\boldmath $A_x$} for some vertices in {\boldmath $V(H_{m+1})$} if {\boldmath $s=m$}.}
We are almost ready to apply Lemma~\ref{lem:technical_bandwidth} if $s< m$. However, if $s=m$, then we need to make sure that we exhaust all colors in $B''$. For this, we will pick an $3$-independent matching $M$ in $H_{m+1}$. We plan to manually embed those edges in $M$ and color them using the colors in $B''$. This will put some restrictions on how the neighbors of the vertices in $V(M)$ should be embedded. To control this, we define the target sets $A_x$ for more vertices $x$.
If $s< m$, then we take $\mathcal{I}=\{ [t_{s}+4k^3+1,t_{s}+4k^3+k-1] \}$ consisting of just one interval. 

In order to define $\mathcal{I}$ for the case of $s=m$, we assume $s=m$ throughout Step 2-2.
Let $t' = t-t_m$, then \eqref{eq: choice of Hsss} implies 

\begin{equation} \label{comparison: t'}
t' > \frac{(\beta_2 -2k\D\g) n}{4\Delta k\alpha n} > \beta_2^2 \alpha^{-1}.
\end{equation}
We enumerate all edges in $R^{m+1}$ into $g^*_1,\dots, g^*_p$. By Lemma~\ref{lem: clique exist}, we can choose cliques $Q^*_1,\dots, Q^*_p$ such that $g^*_i$ belongs to $Q^*_i$ for each $i\in [p]$. 
We know that $r_{m+1}\leq  e(R^{m+1}) = p \leq \binom{r_{m+1}}{2}$.
Consider intervals $\mathcal{I} =\{ I_0,\dots, I_{p} \}$ of $[t_{m}+4k^3+1,t]$ such that the following hold.
\begin{enumerate}
    \item[(a)] $|I_i| = \beta_1 t'$ for all $i\in [p]\cup \{0\}$ and $t_{m}+4k^3+1 \in I_0$.
    \item[(b)] $\mathcal{I}$ is $2\beta_1 t'$-fragmented.
    \item[(c)] $\mathcal{I}$ is $\beta_1^{1/4} t'$-initial. In other words, $I_i\subseteq [t_{m}+4k^3+1, t_{m}+4k^3+ \beta_1^{1/4} t']$ for all $i\in [p]$.
\end{enumerate}
Indeed, we can simply take $I_{i} = [t_m+4k^3+ (3i-3)\beta_1 t'+1, t_m+4k^3+ (3i-2)\beta_1 t']$ for $i\in [p]$. Then $(3p-1)\beta_1 t' \leq \beta_{1}^{1/4} t'$ as $0< \beta_1 \ll 1/r'_1, 1/k < 1$ and $r_{m+1} \leq r'_1$.
Let $I'_i$  be the subinterval of $I_i$ obtained by deleting the first $k$ numbers and the last $k$ numbers.

\begin{claim}\label{cl: matching M}
There exists a matching $E\subseteq \binom{V^{m+1}}{2}$ and a bijection $\lambda'':E\rightarrow B''$ satisfying the following.
    \begin{enumerate}
        \item[(a)] For each $e\in E$, we have $e\in G_{\lambda''(e)}$. 
        \item[(b)] For each $e=uv\in E$, there is $h(e)\in [p]$ such that $u\in V^{m+1}_{j}$ and $v\in V^{m+1}_{j'}$ for $jj'=g^*_{h(e)}$. Moreover, for any $j''\in V(Q^*_{h(e)}) \setminus \{j,j'\}$, we have 
        $$|N_{G^{m+1}}(\{u,v\})\cap V^{m+1}_{j''}| \geq \frac{d^2}{3} |V^{m+1}_{j''}|.$$
        Moreover, 
        $|N_{G^{m+1}}(u)\cap V^{m+1}_{j'}|\geq \frac{d^2}{3} |V^{m+1}_{j'}|$ and $|N_{G^{m+1}}(v)\cap V^{m+1}_{j}|\geq \frac{d^2}{3} |V^{m+1}_{j}|$  also hold.
    \end{enumerate}
\end{claim}
\begin{poc}
    We can find $E$ by choosing an edge $e\in G_i$ for each $i\in B''$ one by one.
    Assume we have a matching $E'$ with $\lambda''(E')\subseteq B''$ and $i\in B''\setminus \lambda''(E')$.
    Since $\delta(G_i[V^{m+1}]) \geq \left(1-\frac{1}{k} + \frac{\eps}{4}\right) n_{m+1}$, the graph $G_i[V^{m+1}]$ contains at least $\left(1-\frac{1}{k}+\frac{\eps}{4}\right)\binom{n_{m+1}}{2}$ edges.
    
    On the other hand, we have $\delta(R^{m+1})\geq \left(1-\frac{1}{k}+\frac{\eps}{4}\right)r_{m+1}$, hence the number of pairs $uv\in \binom{V^{m+1}}{2}$ which do not lie in $(V^{m+1}_j,V^{m+1}_{j'})$ for some $jj'\in E(R^{m+1})$ is at most 
    $$\left(\frac{1}{k}-\frac{\eps}{4}\right) \frac{r_{m+1}^2}{2} \frac{\left(1+\eps^{1/100}_0\right)^2 n_{m+1}^2}{r_{m+1}^2}  + r_{m+1} \left(\frac{n_{m+1}}{r_{m+1}}\right)^2\leq \left(\frac{1}{k}-\frac{\eps}{8}\right)\binom{n_{m+1}}{2}.$$
    Since $k\geq 2$, we have 
    $$\left(1-\frac{1}{k}+\frac{\eps}{4}\right)\binom{n_{m+1}}{2}  -\left(\frac{1}{k}-\frac{\eps}{8}\right)\binom{n_{m+1}}{2} \geq \frac{\eps}{4} \binom{n_{m+1}}{2}$$
    pairs $uv$ which lie in $G_i$ as well as lie between $V_{j}^{m+1}$ and $V^{m+1}_{j'}$ for some $jj'\in E(R)$. By pigeonhole principle, there exists $jj'\in E(R)$ such that 
    $G_i[V^{m+1}_j, V^{m+1}_{j'}]$ contains at least $\frac{\eps}{4} |V^{m+1}_j| |V^{m+1}_{j'}|$ edges. As $jj'\in E(R^{m+1})$, we can assign $h(e)$ so that $jj'= g^*_{h(e)}$.

    Recall that $Q^*_{h(e)}$ is a clique containing $jj'$.
    So Lemma~\ref{lem: typical pairs} yields that at least $(1- \eps_0^{1/200})|V^{m+1}_j| |V^{m+1}_{j'}|$ pairs in between $V^{m+1}_j$ and  $V^{m+1}_{j'}$ satisfies the condition (b). Since $\gamma \ll 1/r'_1 < 1/r_{m+1}$  and  $|\lambda''(E')|\leq |B''|\leq \gamma^{1/2} n < \eps^2 n/r^2_{m+1}$, we have 
    $$\left(1- \eps_0^{1/200}\right)|V^{m+1}_j| |V^{m+1}_{j'}| + \frac{\eps}{4} |V^{m+1}_j| |V^{m+1}_{j'}| \geq \left(1+\frac{\eps}{8}\right)|V^{m+1}_j| |V^{m+1}_{j'}| > |V^{m+1}_j| |V^{m+1}_{j'}| + 2|\lambda''(E')|n.$$
    Thus, we can choose a pair $e=uv\in G_i[V^{m+1}_j, V^{m+1}_{j'}]$ which does not share any vertex with any pairs in $E'$. By repeating this for all colors in $B''$, we obtain the desired function matching $E$ and a bijection $\lambda'': E\rightarrow B''$.
\end{poc}

Take a matching $E$ and a coloring $\lambda''$ as in the above claim.
For each $j\in [p]$, let ${E_j = \{e\in E: h(e)=j\}}$ and $B''_j = \lambda''(E_j).$
Let $m_j = |E_j| = |B''_j|$, then we have $m_j \leq \sum_{j\in [p]} m_j \leq |B'| \leq \gamma^{1/2} n$ by \eqref{eq: property B'}.

For each graph $H[W_{I_j}]$, we next apply Lemma~\ref{lem: matching M} to obtain a $3$-independent matching of size $m_j$ lying in $I_j'$. Since $H$ has no isolated vertices and $\mathcal{W}$ is a bandwidth ordering, the vertices in $W_{j'}$ are isolated in $H[W_{I'_j}]$ only if $j'$ is either one of the $k-1$ smallest numbers in $I'_j$ or $k-1$ largest numbers in $I'_j$. Hence, $H[W_{I'_j}]$ contains at least $\frac{1}{2}(\beta_1 t'-2k)\alpha n$ edges.
By Pigeonhole principle, we can choose two colors $a_j,b_j \in [k]$ so that at least $\frac{1}{2k^2}(\beta_1 t'-2k)\alpha n$ edges are between the vertex sets $W_{I'_j}\cap Y_{a_j}$ and $W_{I'_j}\cap Y_{b_j}$ some $a_j\neq b_j\in [k]$. In other words, those edges are between the vertices of color $a_j$ and $b_j$ under the proper coloring $c$.
Thus, Lemma~\ref{lem: matching M} yields a $3$-independent matching of size 
$$\frac{(\beta_1 t' -2k) \alpha n}{2 k^2 \Delta^3} \geq \beta_2^3 n > \gamma^{1/2} n > m_j,$$
where the first inequality follows from $t'> \beta_2^2 \alpha^{-1}$ (see, \eqref{comparison: t'}).
Hence, we find a $3$-independent matching $M_j$ of $H[W_{I_j}]$ with $|M_j|= m_j$ and $V(M_j)\subseteq W_{I'_j}$ for each $j\in [p]$, such that the edges in $M_j$ is between $W_{I'_j}\cap Y_{a_j}$ and $W_{I'_j}\cap Y_{b_j}$.

Recall that $g^*_j$ is an edge in $Q^*_j$. Now, for each $j\in [p]$, we turn each clique $Q^*_j$ into an ordered clique by taking an ordering $(q(1),\dots, q(k))$ of it so that it satisfies the following.
\begin{align}\label{eq: g* ajbj aligned}
q(a_j) q(b_j)= g^*_j.    
\end{align}
We now take arbitrary bijections between $M_j$ and $E_j$ for each $j \in [p]$. This yields a map from the vertex set $V(M_j)$ to $V(E_j)$ in such a way that each vertex in $M_j\cap Y_{a_j}$ is mapped to $V_{q(a_j)}$ and each vertex in $M_j\cap Y_{b_j}$ is mapped to $V_{q(a_j)}$.

Doing this for all $j\in [p]$ yields a bijection $\varphi^*$ between the vertex set $V(M)$ of a $3$-independent matching $M= \bigcup_{j\in [p]} M_j$ and the vertex set of $E$, and the coloring $\lambda''$ on the image $\varphi^*(M)=E$ is a bijection. As $M$ is an induced matching, $\varphi^*$  does not map any other edges of $H$.
Furthermore, as the edges in $M$ lie in $\bigcup_{j\in [p]} W_{I'_j}$ and $I'_j$ is obtained from $I_j$ by deleting first and last $k$ sets $W_i$, by the bandwidth of $H$, the vertices outside $\bigcup_{i\in [p]} I_j$ does not have any neighbor in $V(M)$.
Now, for each $j\in [p]$ and a vertex $x\in I_j$, we have to designate a target set where this vertex $x$ can be embedded.  For each $\ell\in [k]$ and $x\in I_j\cap Y_\ell$, we take 
$$A_x = N_{G^{m+1}}(\varphi^*(N_H(x)\cap V(M)))\cap  V^{m+1}_{q(\ell)}.$$
Since $M$ is $3$-independent, $N_H(x)\cap V(M)$ is a subset of some $e={\varphi^*}^{-1}(e') \in M$ for some $e'\in E$, 
Claim~\ref{cl: matching M}-(b) and \eqref{eq: g* ajbj aligned} ensure that for all $x\in I_j$, we have
\begin{align}\label{eq: Ax in I}
    |A_x|\geq \frac{d^2}{3}|V^{m+1}_{q(\ell)}| > d' |V^{m+1}_{q(\ell)}|.
\end{align}
Among the vertices in $W_{I_0}$, for all vertices $x$ from the first $k-1$ sets, we already defined $A_x$. For the rest of $x\in W_{I_0}$, if $x\in W_j$ and $j-t_m = \ell ~({\rm mod}~k)$ for some $\ell\in [k]$, then we take $A_x = V^{m+1}_{\ell}$. We now have defined our collection $\mathcal{I}$ of $\beta_1^{1/4} t'$-initial $2\beta_1 t'$-fragmented intervals and our target sets $A_x$ for the case of $s=m$. 
\newline

\noindent {\bf Step 2-3. Embedding {\boldmath $H^{\rm band}$} by Lemma~\ref{lem:technical_bandwidth}.}
Let $G'= G^{s+1}[(V^{s+1}\setminus V'')\cup U^{s+1}_{[k-1]}]$. In particular, if $s=m$, then $U^{s+1}$ is the emptyset and $G'=G^{m+1}[V^{m+1}\setminus V'']$.
Recall that we have defined a collection $\mathcal{I}$ of $\beta_1^{1/4} t'$-initial $2\beta_1 t'$-fragmented intervals before, which is just $\mathcal{I} = \{I_0\}$ with $I_0 =[t_{s}+4k^3+1, t_{s}+4k^3+k-1]$ if $s<m$.
Now we apply Lemma~\ref{lem:technical_bandwidth} to embed $H^{\rm band}$ into $G^{s+1}[V^{s+1}\setminus V'']$ with the following parameters.
\begin{center}
\begin{tabular}{ |c|c|c|c|c|c|c|c| }
\hline 
 parameters/objects & $|V(H^{\rm band})|$ & $\eps/8$ & $d/2,d'/2$ & $G'$ & $(2\eps_0)^{1/200}$ & $r_{s+1}$ & $\alpha \frac{n}{|V(H^{\rm band})|}$ \\ 
 \hline
 playing the role of & $n$ & $\eps$ & $d,d'$ & $G'$ & $\eps'$ & $r$ & $\alpha$ \\ \hline \hline
   parameters/objects & $\{W_{t_{s}+4k^3+1},\dots, W_{t_{s+1}}\}$ & $[k]$ & $U_i^{s+1}$ & $\mathcal{I}$ &  $\mathcal{M}^{s+1,2}$ & $2\beta_1$ & $V_i^{s+1}\setminus V''$ \\ 
 \hline
 playing the role of & $\mathcal{W}$ & $K$ & $U_i$ & $\mathcal{I}$ & $\mathcal{M}$ &$ \eta$ & $V_i$  \\
 \hline
\end{tabular}
\end{center}
Recall that the five conditions in \eqref{eq: Us+1 property} and \eqref{eq: min deg for next} hold with the slightly modified parameters after deleting~$V''$. Using these and \eqref{eq: size of Ax}, similarly to before, we have the following. 
\begin{enumerate}[leftmargin=*]
    \item[(i)] $G^{s+1}[V^{s+1}\setminus V'']$ admits $((2\e_0)^{1/200},d/2)$-regularity partition $\left(R^{s+1},(V^{s+1}_1\setminus V'', \ldots, V^{s+1}_{r_{s+1}}\setminus V'')\right)$ and  for all $i\in [r_{s+1}]$, we have $\left|V^{s+1}_i \setminus V''\right| = \left(1\pm (2\e_0)^{1/400}\right)\frac{n_{s+1}}{r_{s+1}}$.
    \item[(ii)] $\D(\M^{s+1,2})\le \D$ and the rank of $\M^{s+1,2}$ is at most $\D$.
    \item[(iii)] For $s<m$, for $x\in W_{t_{s+1}+l}$ with $\ell \in [k-1]$, if $N_x\cap W_{t_{s+1}-k+i} \neq \varnothing$ for some $i\in [k]$, then $(V^{s+1}_i\setminus V'',U^{s+1}_\ell)$ is an $((2\e_0)^{1/200},d/2+)$-regular pair in $G'$.
    \item[(iv)] For each $I_j\in \mathcal{I}$, there exists an ordered $k$-clique $Q_j=(q(1),\dots, q(k)) \in \overrightarrow{K_k}(R^{s+1})$ such that the following holds:  for each $x\in W_{j'}$ with $j'\in I_j$ with $j'=\ell ~(\mathrm{mod}~k)$ for some $\ell\in [k]$, then ${|A_x|\ge \frac{d'}{2}\left|V^{s+1}_{q(\ell)}\setminus V''\right|}$. 
\end{enumerate}
Here, (iv) comes from the definition of $A_x$ by taking the clique $Q^*_j$ for the interval $I_j$.
Thus, we conclude that the four conditions \ref{cond: G regularity}-\ref{cond: target sets} hold with the above parameters and objects. Consequently, Lemma~\ref{lem:technical_bandwidth} yields an embedding $\varphi'':H^{\rm band} \rightarrow G^{s+1}[V^{s+1}\setminus V'']$ such that the following two properties hold.
\begin{enumerate}[leftmargin=*]
    \item[(a)] For each $x\in I_{j}$ for some $j\in [p]$, $\varphi''(x)\in A_x$.
    \item[(b)] For each $j\in [k-1]$ and $x\in W_{t_{s+1}+j}$, we have the following where ${A_x = \bigcap_{y\in N_H(x)} N_{G^{s+1}}(\varphi'(y)) \cap U^{s+1}_{j}}$.
    $$|A_x| \geq d' |V^{s+1}_j|$$
\end{enumerate}
The property (a) ensures that $\varphi_{s+1} = \varphi'_s\cup \varphi''$ maps all edges of $H$ between the vertices in $\bigcup_{i\in [p]} I_j$ and the rest into $G^{s+1}$. 
In the graph $\varphi_{s+1}(H_{[s+1]})$, the edges in $\varphi'_s(H[W_{[t_{s+1}]}])\setminus ( \varphi_s(H[W_{t_{s}+4k^3}])\cup E )$ are not yet colored by $\lambda'_s \cup \lambda''$.

For each edge $e\in \varphi'_s(H[W_{[t_{s+1}]}])\setminus (\varphi_s(H[W_{t_{s}+4k^3}]) \cup E)$, as it is an edge of $G^{s+1}=\mathcal{G}_{B_{s+1}}^{\eps/8}$, there are at least $(\eps/8)|B_{s+1}|$ colors $i\in B_{s+1}$ such that $e\in G_i$. As at most $e(H^{\rm con})\leq 4k^3\Delta \alpha n$ colors from $B_{s+1}$ are used before for $\varphi'_{s}$, each such an edge $e$ still have more than $(\eps/8)|B_{s+1}| - 4k^3\Delta \alpha n > 2\Delta \gamma n > e(H^{\rm band})$ many colors available if $s<m$.

If $s=m$, then at most $e(H^{\rm con})+ |E|\leq 2\gamma^{1/2} n$ colors from $B_{s+1}$ are used. By \eqref{eq: property B'}, we have 
$(\eps/8)|B_{s+1}| - 2\gamma^{1/2} n  \geq 2\beta_2 n > e(H^{\rm band}).$

In either cases, we choose a color $\lambda_{s+1}(e)$ for each such an edge $e$ in a greedy way to obtain a partial rainbow coloring $\lambda_{s+1}$ of $\varphi_{s}(H_{[s+1]})$ which extends $\lambda'_s\cup \lambda''$.
Let $\Lambda_{s+1}$ be the set of colors used by $\lambda_{s+1}$. If $s=m$, this yields a partial rainbow coloing $\lambda_{m+1}$ which uses all the colors in $B'$, so only edges in $\varphi_{s+1}(H_0)$ are not colored, and the remaining colors are $A\cup C'$ for some $C'\subseteq C$.

If $s<m$, then we can also check the properties \ref{eq: cond1 for step 2}--\ref{eq: target sets in step 2} holds for $\varphi_{s+1}$. By construction of this embedding $\varphi_{s+1}$, it satisfies \ref{eq: cond1 for step 2}.
Also \ref{eq: cond2 for step 2} is obvious from the choice of $B_{s+1}$ and $\phi_{s+1}$.
\ref{eq: mindeg GiU} follows from \eqref{eq: Us+1 property}-(d) and \ref{eq: Us regularity} follows from \eqref{eq: Us+1 property}-(c). \ref{eq: target sets in step 2} follows from the above property (b) and the definition of $\M^{s+1,2}$. Hence, we obtain a map $\varphi_{s+1}$ satisfying \ref{eq: cond1 for step 2}--\ref{eq: target sets in step 2}. By repeating this for all $s<m$, we obtain a partial rainbow coloring of $\varphi_{m}(H_{[m]})$ satisfying the properties \ref{eq: cond1 for step 2}--\ref{eq: target sets in step 2}. \newline

\noindent{\bf Step 3. Using color absorber to finish the coloring.}
Since we have a partial rainbow coloring $\lambda_{m+1}$ of $\varphi_{m+1}(H)$ which uses all the colors in $B'$, only edges in $\varphi_{s+1}(H_0)$ are not colored at this point, and the remaining colors are $A\cup C'$ for some $C'\subseteq C$. As there are exactly $h=e(H)$ colors, we know $e(H_0)= |A\cup C'|$.
Finally, applying \eqref{eq: color absorption}-(b) with this choice of $C'$ extends $\lambda_{m+1}$ to $\lambda$, which is a bijection between $\varphi_{m+1}(H)$ and $[h]$ such that $e\in G_{\lambda(e)}$ for all $e\in \varphi_{m+1}(H)$. This provides the desired $\mathcal{G}$-transversal isomorphic to $H$ and finishes the proof of Theorem~\ref{main}.

\printbibliography

\end{document}